\newcount\notenumber

\def\note{\advance\notenumber by 1
\footnote{$^{(\the\notenumber)}$}}

\font\tenmsb=msbm10 \font\sevenmsb=msbm7 \font\fivemsb=msbm5
\newfam\msbfam
\textfont\msbfam=\tenmsb \scriptfont\msbfam=\sevenmsb
\scriptscriptfont\msbfam=\fivemsb \edef\msbfam{\ifcase\msbfam 0\or
1\or 2\or 3\or 4\or 5\or 6\or 7\or 8\or 9\or A\or B\or C\or D\or
E\or F\fi}

\ifx\sc\undefined
    \font\sc=cmcsc10 
\fi

\mathchardef\nmid="3\msbfam2D
\def\Bbb#1{{\fam\msbfam\relax#1}}
\def\N{{\Bbb N}}
\def\Q{{\Bbb Q}}

\def\Z{{\Bbb Z}}
\def\G{{\Bbb G}}

\def\Gal{\mathop{\rm Gal}\nolimits}

\def\ord{{\rm ord}}
\def\P{{\Bbb P}}

\def\E{{\cal E}}
\def\H{{\cal H}}

\def\T{{\cal T}}

\def\1{{\bf C}}
\def\A{{\Bbb A}}

\def\F{{\Bbb F}}
\def\O{{\cal O}}

\def\CVD{{\hfill\hfil{\lower 2pt\hbox{\vrule\vbox to 7pt
{\hrule width 6pt\vfill\hrule}\vrule}}}\par}


\font\title=cmr10 scaled 1200

\baselineskip=15pt

\hsize = 15truecm \vsize = 22truecm

\hoffset = 0.4truecm \voffset = 0.7truecm

\centerline{\title  Hilbert Irreducibility above algebraic groups}\bigskip

\centerline{Umberto Zannier}\bigskip\bigskip

\noindent{\sc Abstract}. 
 This paper concerns Hilbert Irreducibility  for covers of   algebraic groups, with results which appear to be difficult to treat by existing techniques.  The  present method    works by first studying irreducibility above `torsion' specializations (e.g.,  over  cyclotomic extensions)  and then descending the field (by Chebotarev Theorem).  Among the results we offer an irreducibility theorem for the fibers, above a cyclic dense subgroup,  of a cover of $\G_{\rm m}^n$ (Thm. 1) and of a power $E^n$ of an elliptic curve without CM (Thm. 2); this  had not been treated before for $n>1$. As a further application,  in the function field context, we  obtain a kind of Bertini's theorem for algebraic subgroups of $\G_{\rm m}^n$ in place of linear spaces (Thm. 3). Along the way we shall prove other results, as a general 
  lifting  theorem above tori (Thm. 3.1). 
\bigskip

\centerline{\sc \S 1. Introduction.} \medskip

This paper is in the context of   the  Hilbert Irreducibility Theorem (HIT in the sequel); we offer some  results on the   lifting of rational points above  algebraic groups, which appear to be difficult to treat by existing techniques.

In the paper by `cover'  we mean a dominant rational map $\pi:Y\to X$ of finite degree between irreducible varieties. (We remark at once that  by shrinking $Y$ and $X$ to Zariski-open subsets of them, we may actually assume for most purposes that $\pi$ is a morphism, and even that it is {\'etale}, i.e. that that there are exactly $\deg\pi$ points above every point of $X$.)\smallskip

Consider 
a cover  $\pi: Y\to X$, defined over a number field $k$.  Basic  questions in  Diophantine Geometry can be formulated in terms of  the lifting of rational points $x\in X(k)$ to $Y$: when does it happen that $\pi^{-1}(x)$ contains a point in $Y(k)$, or is $k$-irreducible?  The classical HIT states that if $X=\A^n$ then one may find $x\in X(k)$ such that  this last possibility happens (even considering simultaneously finitely many covers of $X$). 

Now, it is of interest, also for the applications, to obtain such  `good' points $x$ in restricted sets of rational points. Situations which are obviously relevant occur when $X$ is an algebraic group, because these  are the fundamental varieties where we are able to generate systematically rational points.\note{The original case  of HIT  is no exception:  $\A^n=\G_{\rm a}^n$ as a variety.}  Here we shall  study  the {\it lifting of points in a Zariski dense cyclic  subgroup $\Omega\subset X(k)$}, for $X$ either a multiplicative torus or a power of an elliptic curve:  {\it under a necessary geometrical condition on $Y$} (see Def. below) {\it we  prove  that   $\pi^{-1}(x)$ is $k$-irreducible for each $x$ in a suitable coset of finite index  in  $\Omega$}.

We stress that,   though the literature is rich  of many versions of HIT (see e.g. [BG, \S 9.6], [FJ, Ch. 11], [Sch, \S 4.4], [Se1, Ch. 9], [Se2, Ch. 3]),   the said  basic situations   do not appear to fall into existing methods. After the work of Faltings, Vojta and others, much is known for {\it subvarieties} of commutative algebraic groups, but for {\it covers} of them the situation is   still unsatisfactory  in dimension $>1$,  even in very  simply-stated  cases.

Some results (dealing basically with linear recurrences) came implicitly from the papers [Z], [FZ], whose conclusions were applied in [C] to a HIT over linear algebraic groups.  In fact, these papers contain primordial {\it ad hoc} versions of   the present method; here, in addition to further results, this is developed in a  more systematic way, also in view of future possible applications.

In particular, here   we consider  the context of abelian varieties (not touched in [Z], [FZ], [C]) from this viewpoint,  an issue which is explicitly mentioned in  the discussion  in Serre's  [Se2], \S 5.4. 

Below we  focus  on the case $X=E^n$. This  case can be better treated   because of   results on the Galois action on torsion points which are better known than otherwise, but  is not   an {\it a priori} limitation of the method.
  

The  method  consists of two main stages and may be very roughly described as follows: \medskip

{\it (A) To use  a suitable (explicit) HIT over a big (cyclotomic) field, of infinite degree over $\Q$.

(B)  To transfer the irreducibility to points over a number field}.   \medskip

\noindent A relevant  issue here is that   a kind of HIT may be  proved directly  over the big field, actually for {\it explicit} specializations at torsion points:  for $\G_{\rm m}^n$ this ingredient has been essentially done in [DZ]  (see Theorem [DZ] below)  hence (A) applies. In the case of abelian varieties step (A) is obtained below in a different way. Then    the {\it transfer} (B) leads to the sought explicit versions of HIT over number fields. This step involves $v$-adic approximation to the torsion points  coming from  step (A), and Chebotarev Theorem (which may be seen as a $0$-dimensional version of HIT). 

 In all of this, it turns out that   the group structure of torsion points is especially relevant for  the location  of  `good' specializations   in   algebraic groups. \medskip

Before stating   some  conclusions,  we introduce   a  simple geometrical condition on the  covers, which shall turn out to be necessary and sufficient for our purposes. 

For $X$ a commutative connected algebraic group, we let $[m]:X\to X$ denote the multiplication map.  By `irreducible' we mean throughout `$\bar k$-irreducible' (supposed for all involved $X,Y$).\medskip
 
 \noindent{\bf Definition:} We say that the cover $\pi:Y\to X$ satisfies condition (PB)  (`pull-back') if for any integer $m>0$ the pull-back $[m]^*Y:=X\times_{[m],\pi} Y$ is irreducible. \medskip
 
\noindent For instance:  if $X=\G_{\rm m}^r$ and $Y:f(x_1,\ldots ,x_r,y)=0$, (PB) means that $f(x_1^m,\ldots ,x_r^m,y)$ is irreducible for all $m>0$.  For $X=\G_{\rm a}^s$, one  instead  finds that (PB) is always  trivially verified.\medskip
 
 \noindent{\bf Remarks.}  (i) Note   that this condition is  unavoidable   for our lifting issues. In fact, suppose that $[m]^*Y$ is reducible, equal to the union  $U\cup V$ of proper closed subsets. Let $\Omega$ be any finitely generated  (dense) subgroup of $X(k)$. By enlarging $k$ to a finite extension, assume that $U,V$ are defined over $k$ and that $\Omega\subset [m](X(k))$. Let $x\in\Omega$ and write $x=[m]x'$ for $x'\in X(k)$. If $\pi(y)=x$, the pair $(x',y)$ is in $[m]^*Y(k)=U(k)\cup V(k)$. For `general'  $x$, this yields a nontrivial splitting of the fiber $\pi^{-1}(x)$ into  two subsets defined over $k$, so $\pi^{-1}(x)$ cannot be $k$-irreducible.

(ii) In   Proposition 2.1   we shall prove in a simple way two  equivalences for this condition.     First, we shall see that it holds if and only if it holds for $m=\deg \pi$ (so it  is a `computable' condition). Secondly, we shall prove that  it holds  if and only if  the map $\pi$ has no nontrival isogeny factors, which shows  the relevance of ramification. (Note that when $\pi$ itself is an isogeny up to birationality, for large  $k$ we have $\Omega\subset \pi(Y(k))$,  so if $\deg\pi>1$ the  irreducibility of $\pi^{-1}(x)$ badly fails.) \smallskip

Let us now give some statements,  starting  with the case $X=\G_{\rm m}^r\times \G_{\rm a}$; we stress that the crux is represented by the $\G_{\rm m}^r$-component (the $\G_{\rm a}$ being included for completeness).\medskip

\noindent{\bf Theorem 1.}   {\it For $i=1,\ldots ,h$, let $\pi_i:Y_i\to X:=\G_{\rm m}^r\times \G_{\rm a}$  be a cover     satisfying (PB). Then, if $\Omega$ is a cyclic Zariski-dense subgroup of $X(k)$, there exists a coset $C$ of finite index in $\Omega$ such that  for all $x\in C$ and for all $i=1,\ldots ,h$  the fiber $\pi_i^{-1}(x)$ is irreducible over $k$.}\medskip

\noindent This result, derived from  the more general Proposition 3.1,  immediately implies a sharp form of the so-called {\it Pisot $d$-th root conjecture} on linear recurrences (proved in   [Z]). It  also implies the more general results  on recurrences  of   [FZ], which   have been   used in [C] to derive an elegant   version of HIT over linear algebraic groups $X$; in this last paper  it is proved in particular ([C], Cor. 7.15) that if  $Y$ is  smooth and  $\pi:Y\to X$ is  finite, then either it  is unramified  or any Zariski-dense semigroup $\Omega\subset X(k)$  contains `good' points; from [C]  (which works also with reducible $Y$) with some work  one may derive a weak version of Theorem 1 in which $C$ is just an infinite set.\note{The paper [C] reduces to  $\G_{\rm m}$ and $\G_{\rm a}$ by  considering   subgroups  generated by a single matrix; the component of the identity in the closure  
of such a subgroup is isomorphic to 
$\G_{\rm m}^r\times \G_{\rm a}^e$,  $e=0,1$.}

In the context of abelian varieties, we have the following analogue for powers of an elliptic curve $E$ without  CM:\medskip

 \noindent{\bf Theorem 2.}  {\it For $i=1,\ldots ,h$, let $\pi_i:Y_i\to E^n$  be a cover   satisfying (PB). Then, if $\Omega$ is a cyclic Zariski-dense subgroup of $E^n(k)$, there exists a coset $C$ of finite index in $\Omega$ such that  for all $x\in C$ and for all $i=1,\ldots ,h$  the fiber $\pi_i^{-1}(x)$ is irreducible over $k$.}
  \medskip

The case $n=1$  follows in stronger `finiteness' form from Faltings' solution of Mordell's conjecture (see  [Se2], \S 5.4),  but for  $n>1$ even  the weaker assertion in which $C$ is just an infinite subset of $\Omega$ had not been treated before.

Theorem 1 looks similar, but can be obtained more rapidly, due to our results for cyclotomic fields  for which we have no counterpart for fields generated by torsion points of abelian varieties (see \S 2). So, Theorem 2 requires additional    arguments, and we treat it separately. 

Theorems 1, 2 concern irreducibility of the fiber  $\pi^{-1}(x)$. Another question is whether a fibre contains rational points. This  easily reduces to the former; an explicit statement is Theorem 4 in \S 4, where (PB) is replaced by the weaker condition of being not  birationally equivalent to an isogeny. \medskip

Our next result is a simple application of the method in the function field context: we offer   a   toric analogue of Bertini Theorem, where algebraic subgroups of $\G_{\rm m}^n$ replace linear subspaces.  For this result we denote by $\kappa$ an algebraically closed field of characteristic zero, and   by $\theta G$ the translate of the algebraic subgroup $G$ by the torsion point $\theta$.
 \medskip
 
 \noindent{\bf Theorem 3.}  {\it  Let  $\pi:Y\to \G_{\rm m}^n$ be a cover defined over $\kappa$, satisfying (PB). Then there is a  finite union $\E$ of  proper connected  algebraic subgroups  of $ \G_{\rm m}^n$ such that if a connected  algebraic subgroup $G$ is not contained in $\E$, then $\pi^{-1}(\theta G)$ is  irreducible (over $\kappa$)  for every torsion point $\theta$.}
\medskip

\noindent The Bertini Theorem may be seen as a version of a similar   statement for  $\G_{\rm a}^n$; a main difference is that in the present case the algebraic subgroups form a `discrete', rather than algebraic,  family, with degrees tending to infinity. We note that also here condition (PB) cannot be omitted. By specialization the theorem may be readily  reduced to the crucial case $\kappa=\overline\Q$.

Here is a  {\bf polynomial version} of the theorem:  {\it Let $f\in \kappa [x_1,\ldots ,x_n,y]$ be of   degree $d>0$ in $y$ and such that $f(x_1^d,\ldots ,x_n^d,y)$ is   irreducible. Then there is a finite union $\H_f$ of proper subgroups of $\Z^n$ such that if $(a_1,\ldots ,a_n)\in \Z^n\setminus \H_f$,  then $f(\theta_1t^{a_1},\ldots ,\theta_nt^{a_n},y)\in \kappa[t,t^{-1},y]$ is  irreducible  for all roots of unity $\theta_1,\ldots ,\theta_n$.}  

\noindent
In particular,  the  {\it Kronecker's substitution} $(x_1,\ldots ,x_n)\mapsto (t,t^m,\ldots ,t^{m^{n-1}})$ preserves the  irreducibility over $\overline\Q$ of a polynomial  $f$ as above, for all integers $m$ large enough in terms of $f$. (We wonder whether it suffices that $m>M_0(\deg f)$.) For   results  in the same direction, but  only over $\Q$, not $\overline\Q$, and with an additional assumption on $f$ (not to be self-inversive in the $x_i$), see [Sch3].
\medskip

 We have stated some applications of the method, but along the way we shall obtain other results, and the organization of this paper is as follows. We shall soon conclude this introduction with a few further examples and remarks on the above theorems. In \S 2 we shall discuss condition (PB) and state (relying on [DZ]) a basic ingredient  (Theorem  2.1), an {\it explicit}  HIT  over cyclotomic fields; an analogue for abelian varieties seems not free of interest  but cannot be obtained with the same methods; in this direction we shall formulate a conjecture, related to the Manin-Mumford conjecture.   In \S 3 we 
 shall prove  a Theorem  3.1 on the lifting of rational points  in $v$-adic neighborhoods of torsion points,   leading  to  Theorem 1. 
 In \S 4 we shall obtain Theorem 2,  which  shall be distinctly more delicate than the toric case.  In \S 5 we shall present a brief deduction  of Theorem 3 from Theorem 2.1 and also a different proof of a (weaker but more laborious) version of this for   $E^n$ in place of $\G_{\rm m}^n$. 
  \medskip

   
 \noindent{\bf Further remarks and examples.}  The case of $\G_{\rm a}$ of Theorem 1 (i.e. $r=0$) reduces to a refined version of the classical HIT, obtained first by Schinzel [Sch2], who  proved in any dimension the existence of whole arithmetical progressions of `good' integer specializations. For $\G_{\rm a}$, the fundamental case occurs in dimension $1$, since the closure of a cyclic subgroup of $\G_{\rm a}^s$ is a line. (Note also that (PB) is always true for covers of $\G_{\rm a}^s$.)
 
 On the contrary, the case of $\G_{\rm m}^r$ does not admit such a reduction to curves, because   a  cyclic group  $\Omega=\xi^\Z$  may well be Zariski-dense in $\G_{\rm m}^r$: this happens when the coordinates $\xi_1,\ldots ,\xi_r$ of $\xi$ are multiplicatively independent. 
 In practice, the situation of Theorem 1 without $\G_{\rm a}$ components leads  to diophantine equations of the shape $f(\xi_1^n,\ldots ,\xi_r^n,y)=0$,   $n\in\Z,y\in k$. For $r=1$ one can use Siegel's theorem or other strong results on curves (see [D] or [CZ]) to prove even finiteness in the ramified case. 
 But for $r>1$   finiteness is known only in special cases. See [CZ], Theorem 2, for  a proof of finiteness  on the so-called {\it dominant root assumption}; this is satisfied `often',  but not generally. 
 
 Here is an amusing {\bf example-problem}: Take $Y:\{y^2=x_1+x_2+1\}$,   and $X:=\G_{\rm m}^2$, and let $\pi(x_1,x_2,y)=(x_1,x_2)$ be  the projection $(x_1,x_2,y)\mapsto (x_1,x_2)$.  Take also  $\Omega=(2+i,2-i)^\Z$ where $i^2=-1$.  I do not know of any method to prove finiteness of solutions $(n,y)\in\Z\times k$ of $y^2=(2+i)^n+(2-i)^n+1$.  Theorem 1 yields a whole   progression of  integers  $n$ such that $(2+i)^n+(2-i)^n+1$ is not a square in $k$. 
 
 It may be worthwile to point out why in these cases we indeed expect only finitely many solutions. Consider more generally  the points of $Y$ with  $S$-units  $x_1,x_2\in \O_{k,S}^*$  and  $y\in\O_{k,S}$. We may embed $X=\G_{\rm m}^2$ in $\P_2$ and take the closure $\overline Y$ of $Y$ in $\P_2\times \P_1$. We are seeking the integral points on $\overline Y$ with respect to the complement $D=\overline Y\setminus Y$, which  is the pull-back by $\pi$ of the three lines $\P_2\setminus \G_{\rm m}^2$;  it turns out that   $K_{\overline Y}+D$ is the class of the  ramification divisor, i.e. the closure of  $x_1+x_2+1=0$. This is `big', so a conjecture of Bombieri-Lang-Vojta (see Cor. 4.2, p. 223 of [L2], or Conj. 14.3.2, p. 483 of [BG]) predicts that the integral points all lie on a curve. Given this, it is then easy to derive  finiteness of our actual solutions from Siegel's theorem. (One can also use probabilistic considerations. Further evidence is provided by the function-field case: see [CZ3], Thm. 1.1.)

 \smallskip
 
As to  the  abelian context of Theorem 2, our proofs  are more laborious and represent one of the main  points of this paper. They involve, in place of cyclotomic fields,  a kind of   weak form of step (A)  for fields generated by torsion points of elliptic curves. In principle, these techniques seem to be   extendable  to   other abelian varieties, provided  the Galois properties of the corresponding torsion points are sufficiently well known. 
 
 We have treated cyclic groups for simplicity, but the method extends to finitely generated ones, with some complication of details but with no conceptual difference with  the cyclic case. \smallskip
  
 As to Theorem 3, the cover  of $\G_{\rm m}^2$ given by  $y^2=1+2x_1+x_2$  shows  that   $\E$ cannot be generally taken  $\{(1,1)\}$. 
By standard arguments one can deduce that, on the same assumptions, if $G\not\subset\E$ the set    $\{c\in\G_{\rm m}^n(\bar k):\pi^{-1}(cG)\  \hbox{\rm is reducible}\}$
  is  closed and proper  in  $\G_{\rm m}^n$. However the irreducibility of $\pi^{-1}(cG)$ for arbitrary $c$ is more delicate and to our knowledge not yet  completely clarified (see [CZ2], \S5, for the case $n=2$).  
  
  I owe to a referee that Theorem 3  is related to a theorem of Kleiman (see [K], Thm. 2), which implies in particular that {$\pi^{-1}(gG)$ is regular, for $g$ in an open dense $U=U_G$}. Now, if we also know that $\pi^{-1}(gG)$ is connected, we can deduce irreducibility (as is done in some proofs of Bertini's theorem). It is also to be noted that  the set $U_G$ depends on $G$, whereas the above  $\E$ does not.

    \smallskip

  All  arguments of this paper are effective (except that  for the application of Serre's theorem on the Galois image we need an effective version of it).
  \smallskip 
  
  As  a related topic  for the interested reader, we mention the paper  [CS]: it   deals with Hilbert Irreducibility in the context of linear algebraic groups,  but with an entirely different viewpoint,    focusing  on the role of the ground field.  
  \medskip

 \noindent{\bf Acknowledgements.} It is a pleasure to thank P. Corvaja for several very helpful discussions and advice and for  some  issues motivating  Theorem 2. I  thank E. Bombieri for pointing out important inaccuracies in a previous version and for several comments on the presentation. I further  thank   M. Fried for  other important  comments and M. Bertolini, D. Bertrand  and P. Parent for their kind clarifications  of some  facts on elliptic Kummer theory and for some references. Finally, I thank anonymous referees for helpful advice which considerably improved the exposition. 
\bigskip

\centerline{\sc \S 2. The condition (PB) and a cyclotomic HIT.} \medskip

In this section we prove some simple properties of the condition (PB), and then  deduce a cyclotomic version of HIT from results in [DZ].  As above, $X$ denotes a commutative  algebraic group. We assume it is defined over a number field $k$, and that it is connected over $\bar k$. \medskip

\noindent{\bf Proposition 2.1.} {\it Let $\pi:Y\to X$ be a cover of degree $d$, defined over $\bar k$. Then it satisfies (PB) if and only if $[d]^*Y$ is irreducible. The map $\pi$ factors as $\lambda\circ \rho$, where $\lambda:Z\to X$ is an isogeny of algebraic groups and $\rho:Y\to Z$ is a rational map satisfying  (PB).}\medskip

\noindent{\it Proof.}  We shall work with varieties and maps defined over $\bar k$. Let us consider a decomposition $[m]^*Y=U_1\cup\ldots \cup U_s$ into irreducible components, where $m\ge 1$. The kernel $T$ of $[m]$ operates by translation on $[m]^*Y$, as  $t:(x,y)\mapsto (x+t,y)$, for $t\in T$. Hence $T$ permutes the $U_i$; let $T_1$ be the stabilizer of $U_1$, and set $X_1:=X/T_1$. Note that $X_1$ is an algebraic group isogenous to $X$, and $[m]$ factors as $\lambda\circ \tau$, where $\tau:X\to X_1$ is the projection. Also, the degree of $[m]^*Y\to Y$ is $\deg[m]$, whereas  the $U_i$ are covers of $Y$ of degree  $|T_1|$ and we have $\deg [m]=s\cdot |T_1|$. Let now $(x,y)\in U_1$; the class $x+T_1$ depends only on $y$, so we have a rational map $\eta: y\to x+T_1=\tau(x)\in X_1$. Hence $\eta(y)=\tau(x)$, so $\pi(y)=[m]x=\lambda\circ\tau(x)=\lambda\circ \eta(y)$.

Note that $\deg\lambda = |T|/|T_1|=s$, and this divides $d=\deg\pi=\deg\eta\cdot \deg\lambda$. Hence the map $[d]$ on $X$ factors through $X_1$, and we may write $[d]=\lambda\circ \tilde\tau$, where $\tilde\tau:X\to X_1$ is another isogeny. 

Now, let $u\in \ker\lambda$ and consider the map $\eta_u= \eta+u$ from $Y$ to $X_1$. We have $\pi=\lambda\circ  \eta_u$ for each $u$, and $[d]^*Y=\cup_{u\in\ker\lambda}  X_1\times_{\tilde\tau,\eta_u}Y$. This is a decomposition into the union of $s$ closed proper subsets, proving that if $s>1$, then already $[d]^*Y$ is reducible. 

Finally,  take a factorization as in the statement with $\deg \rho$ as small as possible. If this cover $\rho$  does not satisfy (PB) then the above argument shows that $\rho$ factors nontrivially as $\tilde \lambda\circ \tilde\eta$, for a rational map $\tilde\eta:Y\to Z_1$ and an isogeny $\tilde\lambda:Z_1\to Z$; we have $\pi=(\lambda\circ \tilde\lambda)\circ\tilde\eta$, contradicting the minimality of $\deg\rho$.\CVD
\medskip


From now on we shall tacitly use  the content of this proposition. We go on to state a HIT over cyclotomic fields,  with explicit specializations at the set $\T_r$ of torsion points of $\G_{\rm m}^r$.
We denote by $k^c$ the field generated over $k$ by all roots of unity. By {\it torsion coset} of $\G_{\rm m}^r$ we mean a translate of an irreducible algebraic subgroup by a torsion point (see [BG, Ch. 3] for the simple theory).\medskip

\noindent{\bf Theorem 2.1.} {\it  Let $\pi:Y\to \G_{\rm m}^r$ be a cover defined over $k^c$ and satisfying (PB). Then there exists a finite union $\E$ of proper  torsion cosets such that if $\zeta\in \T_r\setminus \E$ then $\zeta\in\pi(Y)$ and if $\pi(u)=\zeta$, then $[k^c(u):k^c]=\deg\pi$ and $\pi^{-1}(\zeta)$ is $k^c$-irreducible.}
\medskip

{\bf $\bullet$} The argument  in Remark(i) to  the Definition shows that   (PB) is necessary for the conclusion.\medskip

\noindent{\it Proof.}  The crux of the proof lies in the following result, obtained in [DZ] (see Theorem 1 therein): 
\smallskip

\noindent{\bf Theorem [DZ].} {\it Let $Y$ be a $k^c$-irreducible variety and let $\pi:Y\to \G_{\rm m}^r$ be a cover  defined over $k$.  Suppose that $\pi(Y(k^c))\cap \T_r$ is Zariski-dense in $\G_{\rm m}^r$. Then there exists an isogeny $\rho:\G_{\rm m}^r\to\G_{\rm m}^r$ (over $k^c$) and a birational map $\psi:Y\to\G_{\rm m}^r$ such that $\pi=\rho\circ\psi$.} 
\smallskip

We proceed with the proof of Theorem 2.1, denoting $X:=\G_{\rm m}^r$. We can freely replace the field $k$ with a finite extension and  we suppose to have done this so that the finitely many varieties which appear are irreducible over $\bar k$ and defined  over $k$. Let  $\widehat Y$ be a  Galois closure   of $Y\to X$. Then $\widehat Y$ is an irreducible variety  and a Galois  cover of $X$, with group $G$, say. 
Let $x\in X$; if $x$ lies out of a fixed proper subvariety $W$ of $X$,  
 the fiber $\pi^{-1}(x)$ will have $d:=\deg\pi$ points in $Y$; let  this fiber be $\{y_1,\ldots ,y_d\}$. Then $G$ acts as a transitive permutation group on such  fiber.  Note that the fiber above $x$ in $\widehat Y$ may be thought as a set of orderings of $\{y_1,\ldots ,y_d\}$, precisely the $G$-orbit of one such ordering, say $(y_1,\ldots ,y_d)$. Let now $\zeta\in X\setminus  W$ be a torsion point and let $(u_1,\ldots ,u_d)\in\widehat Y$ be an ordering of the fiber above $\zeta$ in $Y$; since $\widehat Y$ is defined over $k$, the Galois group $H_\zeta$ of $k^c(u_1,\ldots ,u_d)/k^c$ sends this ordering of $\{u_1,\ldots ,u_d\}$ in orderings  which must correspond  to points of  $\widehat Y$. 
In turn, an  ordering corresponds to a point of $\widehat Y$ if and only if it lies in  the $G$-orbit of $(u_1,\ldots ,u_d)$. Therefore $H_\zeta$ acts as a ({\it decomposition}) subgroup of $G$.

For any subgroup $H$ of $G$, let now $Y_H:=\widehat Y/H$, i.e. the space of $H$-orbits of points of $\widehat Y$; this is a variety  whose function field $k^c(Y_H)$ is the fixed field of $H$ in $k^c(\widehat Y)$. Note that we have a natural map $\pi_H:Y_H\to X$ induced by $\pi$. 
Then, by the above, $\zeta$ lifts to a $k^c$-rational point of $Y_{H_\zeta}$ (namely, the image of $(u_1,\ldots ,u_d)$ in $Y_{H_\zeta}$).  Let us now fix $H$ and consider those $\zeta$ with $H_\zeta=H$.

If $H$ acts transitively on $\{u_1,\ldots ,u_d\}$, then the degree over $k^c$ of any point $u_i$ 
in the fiber is $d=\deg\pi$, and we have the conclusion. Hence, let us suppose in the sequel that $H$ is not transitive as a permutation group on $\{u_1,\ldots ,u_d\}$. Note that this implies that the fiber product $Y_{H}\times_{(\pi_H,\pi)} Y$ is $k^c$-reducible. 

By Theorem [DZ] applied to $Y_H,\pi_H$, and since we are assuming that $\zeta$ lifts to a $k^c$-rational point of $Y_H$,  either there is a proper subvariety $W_H$ of $X$ containing all of these  $\zeta$ or there exists an isogeny $\rho=\rho_H:X\to X$ and a birational map $\psi=\psi_H:Y_H\to X$ such that $\pi_H=\rho\circ\psi$.  
But in this second case  the fiber product $X\times_{(\rho_H,\pi)} Y$ would be reducible, like $Y_{H}\times_{(\pi_H,\pi)} Y$. But since any  isogeny is a factor of some multiplication map, we deduce that $Y$ would not satisfy (PB), contrary to the assumptions.
Hence this case cannot occur, and the exceptional torsion points  $\zeta$ are confined in the proper subvariety $(\bigcup_HW_H)\cup W$. But the Zariski-closure of a set of torsion points in $X$ is anyway a finite union of torsion cosets (see [BG, Thm. 4.2.2, p. 95]). This proves the result, with $\E$ equal to such Zariski-closure.\CVD
\medskip

An extension of Theorem [DZ],  and its consequence Theorem 2.1,  to abelian varieties in place of $\G_{\rm m}^r$ would be  desirable.   We explicitly make the following:\medskip
 
 \noindent{\bf Conjecture.} {\it Let $A/k$ be an abelian variety and $T$ be  its set of torsion points, generating the field  $k(T)$   over $k$. Let $\pi:Y\to A$ be a cover, and suppose that $\pi(Y(k(T)))\cap T$ is Zariski-dense in $A$. Then there exist an isogeny $\rho: B\to A$ and a birational map $\psi:Y\to B$   such that $\pi=\rho\circ\psi$.}\medskip

It turns out that the same arguments of the proof of Theorem 1 in [DZ] do not work in this abelian case,  already for elliptic curves.  

As another motivation for this conjecture, we sketch  how it implies the Manin-Mumford conjecture (proved by Raynaud in 1983) that {\it a curve $C$ of genus $g\ge 2$, embedded in its Jacobian $J$, contains only finitely many torsion points of $J$}. To deduce this from the conjecture, let $\pi:C^g\to J$ be the  map $(P_1,\ldots ,P_g)\mapsto P_1+\ldots +P_g$; it is a surjective (ramified) cover of degree $g!$. If $C$ contains infinitely many torsion points, then  $C^g$ has a Zariski-dense set of points defined over $k(T)$, sent to $T$ by $\pi$. So  the assumptions of the conjecture are verified, and then  let us assume its conclusion. Then the birational map $\psi:C^g\to B$ would be a morphism ([BG], Cor. 8.2.22, p. 238) and $\pi$ would be unramified, a contradiction for $g\ge 2$. 
\bigskip

\centerline{\sc \S 3. A lifting theorem and  applications to  HIT for covers of algebraic tori.} \medskip

We now present a lifting theorem, crucial for  Theorem 1. The  proof illustrates the combination  of parts (A), (B) of  the method.  To state this result, we denote by   $|\cdot |_v$    the sup-norm   with respect to a place $v$ and by $\T_r$ the set of torsion points of $\G_{\rm m}^r$:\medskip

\noindent{\bf Theorem  3.1.} {\it Let    $\pi:Y\to \G_{\rm m}^r$ be a cover defined over $k$,  of degree $d:=\deg\pi >1$ and satisfying (PB).  Then there is  a finite union $\E\subset\G_{\rm m}^r$ of proper  torsion cosets with the following property:  if $\zeta\in \T_r\setminus \E$   there exists a set of positive Dirichlet density of   places $w$ of   $k(\zeta)$,  of residual degree $1$ above $\Q$,  such that $\pi(Y(k(\zeta)))$ does not intersect the set $\{x\in k(\zeta)^r: |x-\zeta|_w<1\}$.}\medskip

\noindent{\bf Remarks.}   
(i) Note that the prime $l:=w|_\Q$ splits completely in $\Q(\zeta)$, because its residual degree there is $1$. In particular,  the set $\{x\in (\Q^*)^r: |x-\zeta|_w<1\}$  contains a whole residue class  in $\Z^r/l\Z^r$.  For instance, if $r=1$, if $\zeta$ has order $m$ and if $\xi\in\Q^*$ has order $h$ modulo $l$, for some $a$ coprime to $h$  the set contains the  powers $\xi^{a+hm}$, all $m\in\Z$. Similar examples shall lead to Theorem 1. 

(ii) Inspection shows that given $Y,\pi$, one may calculate:  equations for the set $\E$, roots of unity $\zeta$ and places $w|l$ with the relevant properties.\medskip

\noindent{\it Proof of Theorem 3.1.}  We may apply Theorem 2.1  to $Y,\pi$, so let $\E$ be the finite union of proper torsion cosets mentioned there. There is a proper subvariety $\E'$ of $\G_{\rm m}^r$ such that the  fiber of $\pi$ outside $\E'$ has exactly $d$-elements (even in a projective closure of $Y$). The Zariski-closure of the set of  torsion points in $\E'$ is another finite union of torsion cosets, and  by enlarging $\E$ we may suppose it  is contained in $\E$. Now, for a torsion point $\zeta\not\in \E$  let $u\in Y(\bar k)$ be such that $\pi(u)=\zeta$.  The conclusion of Theorem 2.1 guarantees that $u$ exists and we have $[k^c(u):k^c]=d$. 

In the sequel, we shall tacitly assume that this is the case for the $\zeta$ in question.
Let $H=H_\zeta$ be   the Galois group of the normal closure $K=K_\zeta$ of  $k(\zeta,u)/k(\zeta)$. 
 (Note  that $K$ depends in fact only on $\zeta$, not on $u$ because $[k^c(u):k^c]=d$, and we have $K=k(\zeta,u_1,\ldots ,u_d)$ where $u_i$ are the elements of $\pi^{-1}(\zeta)$.)

It is a well-known simple fact (attributed to Jordan - see [Se3]) that  $H$ cannot be the union of conjugates of a proper subgroup.\note{A subgroup $B$ of a finite group $H$ has  at most $[H:B]$ conjugates, all of which contain the origin. Hence if $B\neq H$  their union contains $<[H:B]\cdot |B|=|H|$ elements.}  Therefore, since $k(\zeta,u)\neq k(\zeta)$, there exists an element $g=g_{\zeta,u}\in H$ such that 
$u^{g\tau}\neq u^\tau$ for
 all $\tau\in H$.

We now apply the theorem of Chebotarev to the normal closure $K'$ of $K/\Q$. There exists an element $\sigma\in \Gamma:=\Gal(K'/\Q)$ which restricts to $g$ on $K$. In particular, $\sigma$ fixes $k(\zeta)$ pointwise.  We obtain the existence of infinitely many places $l$ of $\Q$  (in fact a set of positive  density), unramified in $K'$ and  such that the Frobenius class of $l$ in $\Gamma$ is the class of $\sigma$. Let then $v$ be a place of $K'$ above $l$ with Frob$(v|l)=\sigma$, and denote by $w$ the place of $k(\zeta)$ below $v$.  We let $\{u_1,\ldots ,u_d\}$ be the fiber of $\pi$ above $\zeta$ and we choose $l$ large enough so that $u_1,\ldots ,u_d$ are defined and remain distinct modulo $v$   (recall that they are distinct) and so that $Y,\pi$ have good reduction at $v$.\note{We need just a simple  concept  of  `good reduction', i.e. we suppose that the reduction of $Y,\pi$ is defined and has  still the same degree. Easy inspection of the proof shows that an estimate $l\ge c^m$, where $c=c(Y,\pi)$ and $m$ is the order of $\zeta$, suffices.}

Since $\sigma$ fixes $k(\zeta)$, the residual degree of $w|l$ is $1$. Let $a\in \G_{\rm m}^r(k(\zeta))$ be such that $|\zeta -a|_w<1$ and consider  the fiber of $\pi$ above $a$.  Suppose that there is an element $b\in Y(k(\zeta))$ so that $\pi(b)=a$. We have  $\pi(b)\equiv \zeta\pmod v$;    
hence the reduction of $b$ at $v$ is defined and $b\equiv u_i\pmod v$ for some $i$. In fact, otherwise the fiber above  the reduction of $\zeta$, in a projective closure of $Y$,  would contain more than $d$ elements and the same would be true for the fiber above $\zeta$ (e.g. by Hensel lifting, or simply by good reduction, on taking $l$ large enough so that $\zeta$ does not lie modulo $v$ in the `exceptional' variety $\E'$ mentioned in the opening argument).

Now,  $b^\sigma=b$, whence $u_i^\sigma\equiv u_i\pmod v$, because $\sigma$ fixes $v$. However any $u_i$ is a conjugate over $k(\zeta)$ of $u$, so of the shape $u^\tau$ for a $\tau\in H$. Hence $\sigma$ does not fix any of the $u_i$ and permutes them, so we would have $u_i\equiv u_j\pmod v$ for some $i\neq j$, a contradiction which proves that   $b$ cannot be defined over $k(\zeta)$,  proving the sought conclusion.
\CVD\medskip

There are several variations on Theorem 3.1, and we mention one of them, useful for  our applications:\medskip

\noindent{\bf Refinement.} {\it Under the same assumptions, let $F$ be a  number field, Galois over $k(\zeta)$
and such that $[F:k(\zeta)]$  is not divisible by any prime smaller than $d=\deg \pi$. Then   we may further prescribe arbitrarily the Frobenius class in $\Gal(F/k(\zeta))$  of the relevant places  $w|l$.}\medskip

\noindent{\it Proof.} With respect to the above arguments, it suffices to observe that $F$ and $K$ are linearly disjoint over $k(\zeta)$, because the degree $[K:k(\zeta)]$ divides $d!$ whereas $[F:k(\zeta)]$ is divisible only by primes $>d$. Hence the Galois group $H^*:=$Gal$(FK/k(\zeta))$ is the product Gal$(F/k(\zeta))\times H$ and we may find an automorphism in $H^*$ which restricts in a prescribed  arbitrary way to $F/k(\zeta)$ and  to $K$. Then the proof works as before, prescribing the action of $\sigma$ also on $F$.  \CVD
\medskip

We now proceed to the proof  of Theorem 1;  we argue through two other statements, which yield  further conclusions.\medskip

\noindent {\bf Proposition 3.1.}\  {\it Let  $\pi:Y\to \G_{\rm m}^r\times\G_{\rm a}$ have degree $>1$ and satisfy (PB) and suppose that   the cyclic group $\Omega$ generated by $\omega:=(\xi,\tau):=(\xi_1,\ldots ,\xi_r,\tau)\in \G_{\rm m}(k)^r\times \G_{\rm a}(k)$ is Zariski-dense in $\G_{\rm m}^r\times\G_{\rm a}$. Then for all large primes $p$ there exist  infinitely many primes $l\equiv 1\pmod p$  such that, for  a class $t_0$ mod $l$ and any  integers $b_1,\ldots ,b_r$ coprime to $p$ and any   $b\in t_0$, $(\xi_1^{b_1{(l-1)\over p}},\ldots ,\xi_r^{b_r{(l-1)\over p}}, b\tau)$ does not lie in $\pi(Y(k))$.

In particular,    there exists a  coset of $\Omega/[l(l-1)]\Omega$  disjoint   from  $\pi(Y(k))$.}\medskip

\noindent{\it Proof.}   In order to apply Theorem 3.1,   we view $\pi$ as a cover of $\G_{\rm m}^{r+1}$; however this might not satisfy (PB);  a simple remedy is to define a modified map by replacing $\pi$ with $\pi':=\lambda\cdot \pi$, for a   fixed  $\lambda\in\G_{\rm m}^r(k)\times\G_{\rm a}(k)$ (where the dot  refers to the group law), supposing that   $\pi'$ is not ramified above the whole $\G_{\rm m}^r\times \{0\}$; all  but finitely many  choices of the last coordinate of $\lambda$   shall do. To check (PB) for $\pi'$, by factoring $[m]$ on $\G_{\rm m}^{r+1}$ through $[m]$ on the first $r$ components and then $[m]$ on the last one, and using that $\pi$ satisfies (PB), we are reduced to show that {\it if a cover of $\G_{\rm m}^{r+1}$ is not ramified above the whole $\G_{\rm m}^r\times \{0\}$, then pull-back by $[m]$ on the last coordinate leaves it irreducible}. Now, if not, then the cover map  would factor through a nontrivial isogeny on the last $\G_{\rm m}$-coordinate (proof of Prop. 2.1). However this isogeny would be a multiplication by a divisor of $m$, and thus would be totally ramified above $\G_{\rm m}^r\times \{0\}$, leading to a contradiction.

Note now that  the denseness of $\Omega$  amounts to  the  $\xi_i$ being  multiplicatively independent elements of $k$ and $\tau\neq 0$.  We   choose $\lambda$ of the shape $(1,\ldots ,1, a\tau)$, with an integer $a$, large enough so that the previous argument applies; hence  we   assume that $\pi':Y\to\G_{\rm m}^{r+1}$ satisfies the assumptions of Theorem 3.1. For a large prime $p$, let us choose a torsion point  $\zeta\in\G_{\rm m}^{r+1}$   of exact order $p$,  satisfying  the corresponding conclusion: note that for large $p$ we may choose it out of the proper  subset $\E$ (relative to $\pi'$). 

By the hypothesis on $\omega$, we have that    for all  large enough  $p$, $\tau$ is a unit at each place above  $p$ and  the coordinates $\xi_1,\ldots,\xi_r$ have multiplicatively independent classes   in $k(\zeta)^*/(k(\zeta)^*)^p$. This independence is not difficult to check: if a product $\xi_1^{a_1}\cdots \xi_r^{a_r}$ is nontrivially a $p$-th power in $k(\zeta)^*$ then it is a $p$-th power in $k^*$ (for $[k(\zeta):k]$ divides $p-1$). For a given $\epsilon >0$, find now with a well-known Dirichlet Lemma an integer $q=q(\epsilon)$ so that the $qa_i$ have not all zero residues $b_i$ mod $p$ satisfying $|b_i|<\epsilon p$; then $\xi_1^{b_1}\cdots \xi_r^{b_r}$ is also a $p$-th power, say $\eta^p$, $\eta\in k^*$, but has height $<n\epsilon p\max h(\xi_i)$, so $h(\eta)<n\epsilon \max h(\xi_i)$. For large enough $p$ one can take an arbitrarily 
small $\epsilon$, which eventually  forces $\eta$ and $\xi_1^{b_1}\cdots \xi_r^{b_r}$  to be   roots of $1$,  contrary to the independence assumption. (See also [Z], Lemma 2.)

Now, the Refinement  applies to $F=k(\zeta,\xi^{1/p})$. Note that, by multiplicative independence modulo $p$-th powers, Kummer Theory shows that $F/k(\zeta)$  is Galois, abelian  of  degree $p^r$.  We thus may find infinitely many primes $l$ and extensions $w$ of $l$ to $k(\zeta)$ such that:\smallskip

(i) The prime $l$ splits completely in $k(\zeta)$.

(ii)  The image  $\pi'(Y(k(\zeta)))$ does not intersect the set $\{x\in k(\zeta)^{r+1}: |x-\zeta|_w<1\}$.

(iii)  The Frobenius of $w$ in $F/k(\zeta)$ equals a prescribed element of $\Gal (F/k(\zeta))$. \smallskip

Now, this Frobenius  is     an  automorphism $g$ fixing $k(\zeta)$ and such that $g(\xi_i^{1/p})=\theta^{h_i} \xi_i^{1/p}$, for some integers $h_i$, where $\theta$ is a primitive $p$-th root of unity; by multiplicative independence modulo $p$-th powers, Kummer Theory again shows that all choices of $h_i$ are possible; if $\zeta=(\theta^{c_1},\ldots ,\theta^{c_r},\theta^c)$ and if $a_i$ are integers coprime to $p$,  we choose $h_i=a_ic_i$. Now, by  (i)  we have $\xi^{l/p}\equiv g(\xi^{1/p)})\pmod v$, where $v$ is a place of $F$ above $w$ with Frobenius $g$, 
so 
by our choice we have   $\xi_i^{b_i{(l-1)\over p}}\equiv \theta^{c_i}\pmod v$, where $b_i$ is any inverse to $a_i$ modulo $p$. Hence, this congruence holds for the place $w$ of $k(\zeta)$ below $v$. 
Also, for large $l$ both $\tau$ and $\zeta$ reduce to $\F_l^*$ modulo $w$, so for $b$ in a whole progression $t_0+\Z l$   we may prescribe that $(b+a)\tau\equiv \theta^c\pmod w$.

 Hence, we have  $(\xi_1^{b_1{(l-1)\over p}},\ldots ,\xi_r^{b_r{(l-1)\over p}},(b+a)\tau)\equiv \zeta\pmod w$ so by (ii) we conclude that $(\xi_1^{b_1{(l-1)\over p}},\ldots ,\xi_r^{b_r{(l-1)\over p}},(b+a)\tau)$ does not lie in $\pi'(Y(k(\zeta)))$, i.e. that $(\xi_1^{b_1{(l-1)\over p}},\ldots ,\xi_r^{b_r{(l-1)\over p}},b\tau)$ does not lie in $\pi(Y(k(\zeta))$.  This proves the first part.

 We now let $b_0$ be an integer coprime to $p$ and find an integer $u_0$ such that $u_0\equiv b_0(l-1)/p\pmod p$ and $u_0(l-1)\equiv pt_0\pmod l$. Putting $b_i=u_0+m_ip$, $b=u_0+ml$   for arbitrary integers $m_1,\ldots ,m_r,m$, we obtain that  $(\xi_1^{u_0{(l-1)\over p}+m_1(l-1)},\ldots ,\xi_r^{u_0{(l-1)\over p}+m_r(l-1)},(u_0{l-1\over p}+ml)\tau)$ does not lie in $\pi(Y(k(\zeta))$.  In turn, for $m_i=nl$, $m=n(l-1)$, we conclude that 
 $[u_0{(l-1)\over p}+l(l-1)\Z]\omega$ is disjoint from $\pi(Y(k(\zeta)))$, as required. \CVD
 \medskip

\noindent{\bf Corollary.}   {\it Suppose that for $i=1,\ldots ,h$, $\pi_i:Y_i\to \G_{\rm m}^r\times \G_{\rm a}$ is a cover of     degree $>1$ satisfying (PB) and suppose that   the cyclic group $\Omega\subset \G_{\rm m}(k)^r\times \G_{\rm a}(k)$   is Zariski-dense in $\G_{\rm m}^r\times\G_{\rm a}$.  Then there exists a coset  $C$ of finite index in $\Omega$ and disjoint from $\cup_{i=1}^h\pi_i(Y_i(k))$.}\smallskip

\noindent{\it Proof.} We   argue  by induction on $h$, for $h=0$ the assertion being empty. Suppose it proved up to $h-1$, and let $[a+q\Z]\omega$ be the corresponding coset  involving the first $h-1$ maps, where $\omega$ generates $\Omega$.  
We now apply the last assertion of Proposition 3.1, with $Y:=Y_r$,  $\pi:=(-a\omega)\cdot \pi_r$ (this still verifies (PB)) and $[q]\omega$ in place of $\omega$; we  obtain $a',q'$ such that for $m\in\Z$,   the point $[a+(a'+mq')q]\omega$ does not lift under $\pi_r$ to $Y_r(k)$. This completes the induction, with $C=[a+a'q+qq')]\Omega$. \CVD
\medskip

\noindent{\it Proof of Theorem 1.} Note that this Corollary is in fact a weak version  of Theorem 1; we are going to use a rather standard method for the converse deduction.

A first point is to ensure (PB) for the Galois closures of our covers;  let us drop for a moment the index $i$ and let us  consider a Galois closure $\widehat Y$ of  $\pi:Y\to X:=\G_{\rm m}^r\times \G_{\rm a}$. (In all of this we can freely enlarge $k$ to a finite extension.) Now let us take the pullbacks $[B]^*Y,[B]^*\widehat Y$ by a multiplication map $[B]$ on $X$, where $B$ is   divisible by $[\widehat Y:X]$. Now, $[B]^*Y$ is irreducible by assumption, but $[B]^*\widehat Y$ may become reducible, and let $V$ be a component,  noting it is   a   Galois closure of $[B]^*Y$ over $X$. The natural map $\pi:\widehat Y\to X$ has degree dividing $B$, hence by the last assertion of Proposition 2.1, applied to $\widehat Y\to X$, we deduce that $V$ satisfies (PB). 

Letting now $\omega$ be a generator of $\Omega$, let us  replace  $Y$ by $[B]^*Y$ and $\omega$ by $\omega^*$, where $[B]\omega^*=\omega$, extending $k$ so that every point in $[B]^{-1}(\omega)$ is defined over $k$; since $\pi$ lifts to a map also denoted $\pi$ from $[B]^*Y$ to $X$, of the same degree,   we are then reduced to prove the assertion with $[B]^*Y$ in place of $Y$ and $\omega^*$ in place of $\omega$. (Note that we can choose a single $B$ which works for all the original covers $Y_i$.) 
Since $V$ is related to $[B]^*Y$ as $\widehat Y$ is related to $Y$, we conclude that we may work under the assumption  that   $\widehat Y_i$ satisfies (PB), for $i=1,\ldots ,h$, as we suppose from now on.  

Now,  dropping again the index, suppose that  a point $[n]\omega$ lifts to a point on $Y$ of degree $<\deg\pi$ over $k$. Then (as in the proof of Theorem 2.1) $[n]\omega$ lifts to a rational point of some $Y_H:=\widehat Y/H$, intermediate inside  $\widehat Y\to X$, where $H$ is some  subgroup of $\Gal(\widehat Y/X)$, intransitive on the fiber on $Y$ of a generic point of $X$; this intransitivity ensures that the degree $[Y_H:X]$ is $>1$. But $\widehat Y/X$ satisfies (PB), and since (PB) clearly transfers to any intermediate cover, we deduce that $Y_H/X$ satisfies (PB) as well. Hence, finally, it suffices to apply the last Corollary  to the larger number of varieties $Y_{i,H}$ so obtained, to deduce that  for a whole arithmetical progression of  $n$ this does not happen. \CVD 
\medskip

 \noindent{\bf Remarks.}  The   theorem of Chebotarev in Theorem 3.1 allows  to {\it transfer} the information  from a torsion point $\zeta$, over the big  field $k^c$, to a point $x$ over $k$, near to $\zeta$ with respect to a suitable place. 
 
   Inspection shows that  Chebotarev's theorem implicitly appears  in    [Sch2], and   was  also developed independently  by Fried [Fr]  for  function fields over finite fields.  
      \bigskip

 \centerline{\sc \S 4. Proof of  the elliptic  HIT  Theorem 2.}\medskip
 
 In this section we shall prove Theorem 2. The general principles are analogous to the proof of Theorem 1 just given. However in the elliptic context we miss Theorem [DZ] and its consequence Theorem 2.1. Hence, step  (A) has to be carried out differently. For this we shall now adopt ideas from [Z], which however need several new ingredients for the present situation; fortunately, this still suffices to provide a partial substitute. We pause for a   brief sketch before the details.  
  
 The approach uses the well-known Lang-Weil estimate for points of varieties over finite fields (derived from Weil's Riemann-Hypothesis for curves). We recall  this in the following form (see also [Se1], p. 184 or [Se2], p. 30): {\it Let $Z/k$ be an absolutely irreducible variety of dimension $n$. For a  prime $p$, let $v|p$ be a place of $k$ with residue field contained in the finite field $\F_q$. Then, as $p\to\infty$, the number $|Z_v(\F_q)|$ of points of the reduction $Z_v$ of $Z$ satisfies $|Z_v(\F_q)|=q^n+O(q^{n-{1\over 2}})$.} 
 
 Let us now  go back to our  setting of a map $\pi:Y\to \G_{\rm m}^n$, and let $\widehat Y$ be a   Galois closure over $k$, assuming it to be irreducible over $\bar k$. As in a method introduced by Eichler, Fried and by S.D. Cohen (who applied it to HIT), this statement,  applied first to $Z=Y$ and then to $Z=\widehat Y$,  allows to show that the image $\pi(Y(\F_q))$ has $< cq^n+o(q^n)$ elements, for a $c<1$, actually $c= 1-(1/d!)$ works. See [Se1], pp. 184/185 or [Se2], Thm. 3.6.2 for details of this deduction. 
 
 With this in hand,  for large enough $q$ we can choose a point in $X(\F_q)\setminus \pi(Y(\F_q))$, and we lift this to a torsion point $\zeta$ of $X(\bar k)$, through a reduction map modulo a place of $\bar k$  of good reduction,  lying above $p$.  Then, by reduction one may easily check   that  $\zeta$ does not lie in $\pi(Y(k(\zeta)))$.  At this point we have an information similar to (although  weaker than)  the conclusion of Theorem 2.1, and this provides a starting point for the step (B). As to step (B), in principle it  is entirely similar to what is  carried out in Theorem 3.1.   For the present application,   additional difficulties come in when  dealing with several covers simultaneously and  with the Kummer Theory which appears in Proposition 3.1; it is here that we use the special abelian variety   $E^n$, for which  we have Serre's results on the Galois  action [Se4]. However, taking for granted the appropriate results from Galois action and Kummer Theory  the method could work generally.  \medskip

\noindent Let us now go on with the details of the proof of Theorem 2; we shall argue through auxiliary facts.  Let from now on $E$ be an elliptic curve defined over $k$, without CM, and set $X:=E^n$. One begins by  a   step exactly analogue to the reduction of Theorem 1 to the Corollary to Proposition 3.1: precisely the same    argument given for  the proof of Theorem 1 after the Corollary  (which works for general commutative algebraic groups) shows that we only need to prove the following statement:\medskip

\noindent{\bf Proposition  4.1.} {\it Suppose that for $i=1,\ldots ,h$, $\pi_i:Y_i\to X$ is a cover of     degree $>1$, such that a Galois closure $\widehat Y_i/X$ of $Y_i/X$ satisfies (PB). Suppose also  that   the cyclic group $\Omega\subset X(k)$   is Zariski-dense in $X$.  Then there exists a coset  $C$ of finite index in $\Omega$ and disjoint from $\cup_{i=1}^h\pi_i(Y_i(k))$.}\medskip

\noindent{\it Proof of Proposition 4.1.}  The  easy inductive argument given above for the Corollary to Proposition  3.1 allows   us to  reduce  to the case $h=1$.  Thus from now on we drop the index $i$.  Also, by enlarging $k$ to a finite extension if necessary, we suppose in the sequel that all the varieties which appear are $\bar k$-irreducible and defined over $k$ and  we shall indicate by  $c_1,c_2,...$ positive numbers (integers if necessary), depending   only on $k,E,\pi$.  

As usual, $E[m]$ denotes the kernel of $[m]$ on $E$, whereas a  tilde shall denote reduction modulo a place. We shall refer by ST to the already cited Serre's   theorem [Se4] that the Galois group of the field generated over $k$ by all the $E[m]$  has  finite index in $\prod_lGL_2(\Z_l)$. In particular, we may choose  $c_1$ such that if $m$ has  all   prime factors  $\ge c_1$ then $\Gal(k(E[m])/k)\cong GL_2(\Z/m)$. 
\medskip
A first step is  to obtain a certain torsion point $\zeta\in X$  such that any lift to $Y$ has the `correct' degree $\deg \pi$ over $k(\zeta)$. 
   To gain this  irreducibility of $\pi^{-1}(\zeta)$ over $k(\zeta)$,  we  may use  again  the  `standard' argument of Theorem 2.1;  we  may  reduce to show that $\zeta$ has no lift over $k(\zeta)$,   to any  among   the  finitely many $Y_H$, defined  as the quotients $\widehat Y/H$, where $H$ is an intransitive subgroup of $\Gal(\widehat Y/X)$.   Note that these varieties are subcovers of $\widehat Y/X$ and thus satisfy (PB) (because $\widehat Y/X$  does, by the present assumptions).  This motivates us to pause by   proving the following: \medskip
  
  \noindent{\bf Lemma 1.} {\it There exists a torsion point $\zeta=(\zeta_1,\ldots ,\zeta_{n})$ of $E^n$ of order divisible only by primes $>c_1$, such that it does not lift to $Y_H(k(\zeta))$, for any of the finitely many $H$ in question.}\medskip

\noindent{\it Proof of Lemma 1.}    We might  argue as in the opening sketch, by reduction to a finite field $\F_q$, applying the Lang-Weil theorem individually for each of the varieties $Y_H$. But  the problem is that we need a point  $x\in \widetilde E^n(\F_q)$ which is `good'  for all the $\widetilde Y_H$  {\it simultaneously}, i.e. does not lift to any $\widetilde Y_H(\F_q)$, and   such a common point $x$ need not exist over a finite field. (Using several primes, one for each $H$,  does not work since we successively need to lift $x$ to  a torsion point  $\zeta\in E^n(\overline\Q)$.) We overcome this serious obstacle in some steps as follows:  we shall  choose  good points  $x_H\in \widetilde E^n$ (relative to $Y_H$) over a finite field, lift them to torsion points $\zeta_H$ over $\overline\Q$ {\it of a same order, independent of $H$} (which is crucial), and at this stage we conjugate the $\zeta_H$ over $k$ to obtain a same point $\zeta$, simultaneously suitable for all the $Y_H$. \medskip

   With this program in mind, it shall be convenient  (for the moment and for  this task only!) to go to the case of curves, i.e. to $n=1$, by restricting the cover  $Y\to E^n$ above a suitable copy of $E$ inside  $E^n$.    In doing  this we want  to preserve our irreducibility assumptions. Suppose that $n>1$ and that for generic $x\in E$, $Y$ becomes reducible above $\{x\}\times E^{n-1}$. Then the product $Y\times_{E^n} V$ is reducible in $t>1$ components, where $V=C\times E^{n-1}$ and $C$ is a suitable smooth curve with a nonconstant map $C\to E$ of degree $t$. If this map is unramified, then it is a factor of  an isogeny $E\to E$, against  our assumptions (PB) on $Y$. So  the map is ramified, say above $x_0\in E$. But then, since two distinct components of $Y\times_{E^n}V$ merge above any point in  the branch locus of $V\to E^n$, the branch locus of $\pi$ contains $\{x_0\}\times E^{n-1}$. 
   By an automorphism of $E^n$ we may assume this is not the case, so $\pi^{-1}(\{x\}\times E^{n-1})$ is not generically reducible, and so it is  irreducible except for finitely many $x\in E$.    Continuing in this way, we may choose inductively  torsion points $\zeta_1,\ldots ,\zeta_{n-1}\in E$, of large  but fixed   prime orders $l_1,\ldots ,l_{n-1}$, so that $\pi^{-1}(\{(\zeta_1,\ldots ,\zeta_{n-1})\}\times E)$ is irreducible. In fact, we may assume that  the   $\zeta_i$ work simultaneously for all the finitely many $Y_H$ in place of $Y$ and, denoting by $W_H$ the irreducible curve obtained by restricting the cover $Y_H$ above 
   $\{(\zeta_1,\ldots ,\zeta_{n-1})\}\times E$, the same method,  applied to $[B]^*Y_H$ in place of $Y_H$ for a suitable fixed $B$ (the same for all $H$),   shows that  actually   we may prescribe that $W_H$ satisfies (PB).   \medskip

   Let us go on by choosing a sufficiently large prime $p>l_1\cdots l_{n-1}$, which shall be the characteristic of the finite field to work with. We choose $p$  splitting completely in $k_1:=k(E[l_1\cdots l_{n-1}])$,   denoting by $v$ a fixed place of $k_1$ above $p$. This place has residual degree $1$ above $p$, so the reduction $\widetilde E$  of $E$ modulo $v$ is  defined over $\F_p$ and the points in $E[l_1]$  have  reductions in  $\widetilde E(\F_p)$ as well.
  In particular, $l_1$ divides $|\widetilde E(\F_p)|$  so $p$ does not divide $|\widetilde E(\F_p)|$ (which is $<2p$).

  As in the    remarks  at the beginning of this proof,  we shall need points $x_H\in \widetilde E$ having a same order, i.e. independent of $H$.  To keep control on the order and prove the existence of such points   we shall need that $|\widetilde E(\F_p)|$ is not divisible by `high' powers of many `small' primes $l$ (which shall be measured through a certain product $\prod_l(1-l^{-1})$).  Possibly this may be achieved for some $p$, by using Analytic Number Theory, but a delicate quantification would anyway be involved; instead, it is possible, and seems simpler,  to work over a finite field $\F_q$ in place of $\F_p$, $q=p^m$,  keeping $p$ fixed and letting  $m$ vary. \smallskip

   We start with some technical choices, whose motivation shall be clearer later.
   
  Write $\gamma_m=|\widetilde E(\F_q)|$, so $\gamma_m=(1-\alpha^m)(1-\beta^m)$, where $\alpha,\beta$ are the Frobenius eigenvalues (see [Si]); note that since $\gamma_1$ is prime to $p$, $\gamma_m$ is also prime to  $p$ for all $m$ in a suitable arithmetical progression $P$. In fact, let $p'$ be a place of $\Q(\alpha)$ above $p$; since $\alpha\beta=p$ and $\alpha,\beta$ are conjugate over $\Q$,   $p'$ has residual degree $1$. Hence  $\alpha^{p}\equiv\alpha$, $\beta^{p}\equiv\beta\pmod{p'}$, so $\gamma_{m+p-1}\equiv\gamma_m\pmod p$ for $m>0$ and we can take $P=1+(p-1)\Z$.

By a similar argument, using that $\alpha^{l^{2h}}$ and $\beta^{l^{2h}}$ are constant modulo $l'^e$ for a prime $l'$ above $l$ in $\Q(\alpha)$ and $h>h(e)$, we may also assume that,  for all $m\in P$, $ \gamma_m$ is not divisible by a large enough fixed power $l_i^e$ of any $l_i$ or any power $l^e$, $e>c_2$, of any other prime $<c_1$: it suffices to take $P=1+M\Z$ with $M$ divisible by $p-1$ and, for each of the said primes $l$,  by $(l^2-1)l^h$ where $h$ is large enough. We tacitly let $m$ run through  $P$ in what follows. \smallskip

 Thinking of $p$ and $P$ as fixed, we denote by $m_0$ an integer in $P$, to be chosen sufficiently large so to satisfy certain properties that we are going to explain; in doing this we shall use asymptotic formulas which have to be understood as holding for $m_0\to \infty$. We set $q_0=p^{m_0}$ and we  let $N=N(m_0)$ be the largest integer so that 
 $$
 B:=N!\ \le q_0^{1\over 10}=p^{m_0\over 10}<(N+1)!.
 $$
 For fixed $p$ and large $m_0$ we have 
 $$
N\log N\sim  ({\log p\over 10})\cdot m_0,\qquad  N\sim ({\log p\over 10})\cdot {m_0\over \log m_0}.
 $$
 We write $B$ as a product $\prod_{l\le N}l^u$ of prime powers $l^u$, where $l\le N$ and  $u=u_{l,N}\ge N/2l$. \smallskip
 
 Suppose that such a prime-power $l^u$ divides three numbers $\gamma_{\mu_i}$ with $\mu_1<\mu_2<\mu_3$ integers  in $P$ and in a certain interval $J=[m_0,m_1]$. Let $l'$ be a place above $l$ in $\Q(\alpha)$. Then, by recalling $\gamma_m=(1-\alpha^m)(1-\beta^m)$, we find that $\ord_{l'}(1-\alpha^{\mu_i})+\ord_{l'}(1-\beta^{\mu_i})\ge u$ for $i=1,2,3$. By symmetry, we may assume that for two indices $i<j\in\{1,2,3\}$, the maximum order on the left is attained by the $\alpha$-term, so $\ord_{l'}(\alpha^{\mu_j}-\alpha^{\mu_i})\ge u/2$. 
Since  $l\neq p$ (by the above choice of $P$),   $\alpha,\beta,p$ are  coprime to $l$, hence $\ord_{l'}(\alpha^{\mu_j-\mu_i}-1)\ge u/2$. Taking the norm $N^{\Q(\alpha)}_\Q(\alpha^{\mu_j-\mu_i}-1)=\gamma_{\mu_j-\mu_i}$,  we get that $\ord_l\gamma_{\mu_j-\mu_i}\ge u/2$, whence $l^u\le |\gamma_{\mu_j-\mu_i}|^2\le (1+\sqrt p)^{4(\mu_j-\mu_i)}\le   (4p)^{2|J|}$, where $|J|$ denotes the length of $J$.
 
 Therefore  
 $ |J|\ge  u\log l/2\log (4p) \ge  N\log l/4 l\log (4p)\ge \log N/4\log (4p)$, where the last inequality follows from $l\le N$.    Let us then suppose in the sequel that 
 $$
{ \log N\over 8\log (4p)}<|J|=m_1-m_0< {\log N\over 4\log (4p)}.
 $$
  With this choice we have proved that for each of the said prime-powers $l^u$,  {\it the interval $J$ may contain at most two integers $m\in P$ such that $\gamma_m$ is divisible by $l^u$}.\smallskip

 For an $m\in J$ let us now define $\phi_m=\prod_{l^u|\gamma_m}(1-l^{-1})$. We have just  checked  that a  power $l^u$ may divide  at most two integers $\gamma_m$ for $m\in J\cap P$, so a prime $l\le N$ contributes to at most two of the $\phi_m$, for $m\in J\cap P$; hence
 $$
 \prod_{m\in J\cap P}\phi_m\ge\prod_{l\le N}(1-l^{-1})^2\gg  (\log N)^{-2},
 $$
 where the right-hand estimate comes from Mertens's theorem    in  elementary Analytic Number Theory (see [I], Thm. 7, p. 22). 
 
 We  thus obtain $|J\cap P|\log\left(\max_{m\in J\cap P}\phi_m\right) \ge -2\log\log N+O(1)$. Then,  taking into account the choice of $|J|$,  it  immediately  follows  that for fixed $p$, $P$ and for every given $\delta>0$, if $m_0$ is large enough  we have $\max_{m\in J\cap P}\phi_m\ge 1-\delta$.  We choose therefore an $m\in J\cap P$ such that this estimate $\phi_m\ge 1-\delta$ is verified, with a  $\delta$ small enough, in terms of $\deg\pi$, to justify the coming arguments. \medskip

  We now let $W\to E$ denote one of the above covers  $W_H\to E$;   its  degree $[W_H:E]$ is $\le [Y_H:X]\le [Y:X]!$. Note that any fixed cover remains irreducible under reduction modulo a place of large enough norm. Then, since $W_H$ satisfies (PB), we may choose $p$ large enough such that a fixed set   of pull-backs of $W_H$ remains irreducible modulo $v$; hence  by Proposition 2.1 if $p$ is large enough  the reduction of $W_H$ at $v$ satisfies (PB). 
So, denoting $W(B):=[B]^*W$ be the  pullback by the multiplication-map $[B]:E\to E$, the reduction of this   cover  is  irreducible.

  We   denote with a tilde such a reduction modulo $v$, and we choose $m$ in $P$ and in an interval $J$ as above, so that $\phi_m$ is maximal and thus $\ge 1-\delta$. As above we put $q:=p^m$.
  
  We  have $m_0\le m\le m_0+|J|$, so the above displayed inequalities yield
  $$
  m=m_0+O(\log m_0).
  $$

  We now apply the Lang-Weil theorem, as in [Se1], pp. 184/185 or [Se2], Thm. 3.6.2 (see also the above sketch) to the reduction of $W(B)$ and its Galois closure over $E$. More precisely, taking into account that $B$ is varying, we shall  apply the Weil Theorem (Riemann Hypothesis) for curves, in place of Lang-Weil, so to have a uniform control of the error term.\note{Here we could argue differently, without reducing to curves,  by combining bounds of Deligne and Bombieri, as noted in [Se2], \S 3.6. This could be crucial in dealing with  simple abelian varieties in place of $E^n$, but here the present method is simpler.} This error term depends on the genera of the involved curves. Now,  note that the genera of $\widetilde W$ and its Galois closure are  bounded, for $W$ running through the $W_H$ and varying $p$. So, by the Riemann-Hurwitz formula,  the genera of $\widetilde  W(B)$ and its Galois closure are  bounded by $\ll \deg_E[B]= B^2\le q_0^{1/5}$, because any isogeny is unramified and so $W(B)$ is unramified above $W$.

  Let $\Lambda$ be the set of   points in $\widetilde  E(\F_q)$ which do not lift to $\widetilde  W(B)(\F_q)$ (under the natural map $\pi$ on $W(B)$); i.e.,  $\Lambda$ is the set of $x\in \widetilde  E(\F_q)$ such that $[B]x$ does not lift to $\widetilde W(\F_q)$. Then by the Weil Theorem  (applied as in the books  referred to above, as recalled in the opening sketch) we deduce, taking into account our estimate for the genera  that   
  $$
 | \Lambda| \ge c_3q.
  $$
  Here this estimate works for all large $m$, where $c_3$ is a positive constant depending only on $\deg\pi$.

  Now, the group $\widetilde E(\F_q)$ is of the shape $(\Z/a)\oplus (\Z/b)$ for integers $a,b$ with $a|b$ and $ab=\gamma_m$. 
  Let $a_1=a/\gcd(a,B)$, $b_1=b/\gcd(b,B)$. The map $[B]:\widetilde  E(\F_q)\to \widetilde  E(\F_q)$ has kernel $K$ isomorphic to $(\Z/\gcd(a,B))\oplus (\Z/\gcd(b,B))$ and image   isomorphic to $(\Z/ a_1)\oplus (\Z/ b_1)$.  \smallskip

  Note that, by the previous choices, the primes $l\le N$ which may divide $a_1b_1$ are such that the corresponding product $\prod_{l\le N, l|a_1b_1}(1-l^{-1})$ of $1-l^{-1}$ is at least $1-\delta$: in fact if $l|a_1b_1$ then $l$ must divide $ab$ to a power superior to the power $l^u$, with which it divides $B$; so, {\it a fortiori}, $l^u|\gamma_m$ and the above estimate $\phi_m\ge 1-\delta$ applies.\smallskip
  
 Let us now  estimate from below  the product $\prod_{N<l|\gamma_m}(1-l^{-1})$,  of  $1-l^{-1}$ over the primes $l>N$ which divide $\gamma_m$. If there are exactly $h$  such primes, and if $N<p_1<p_2<\ldots <p_h$ are the first $h$ primes greater than $N$, we clearly have 
 $$
 \prod_{N<l|\gamma_m}(1-l^{-1})\ge \prod_{i=1}^h(1-p_i^{-1}).
 $$
  On the other hand, $\prod_{N<l|\gamma_m}l\le \gamma_m$, so $\prod_{i=1}^hp_i\le \gamma_m$, whence $\sum_{i=1}^h\log p_i\le\log \gamma_m\le 2m\log p$. Hence, by Chebyshev's elementary estimates in Prime Number Theory (see [I], Ch. I, \S\S 4,5) we have $p_h\ll \sum_{prime \ l\le p_h}\log l= \sum_{prime\ l\le N}\log l+ \sum_{i=1}^h\log p_i\ll N+2m\log p\ll   m_0\log p$, the implied constants being absolute. Therefore, using again Mertens's theorem  we obtain
  $$
  \prod_{i=1}^h(1-p_i^{-1})\ge {\log N\over \log p_h}(1+o(1))\ge {\log N\over \log m_0}(1+o(1))\ge 1+o(1),
  $$
  where the terms $o(1)$ tend to $0$ as $m_0$ grows to $\infty$ (recall that we are working with a  fixed $p$). Hence, for large enough $m_0$ we find that 
$\prod_{N<l|\gamma_m}(1-l^{-1})$  is also at least  $1-\delta $.

  In conclusion, we have shown in particular that, denoting by $\phi$ the Euler's function,  
  $$
 \phi(b_1)\ge  (1-\delta)^2b_1.
 $$

  Note that since $W(B)$ is a pull-back by the map $[B]$,  the set $\Lambda$ is invariant by addition of $K$, i.e. $\Lambda + K=\Lambda$.
  If every element $(t,u)$ of $\Lambda$ (we refer here to the above direct sum representation)  is  such that $
  [B](t,u)$ has entries $(t_1,u_1)\in \Z/a_1\oplus \Z/b_1$ such that  
  $\gcd(b_1,u_1)>1$, then 
  $$
  |\Lambda|\le |K|(a_1b_1-a_1\phi(b_1))\le (a/a_1)(b/b_1)a_1(1-(1-\delta)^2)b_1\le 2\delta ab.
  $$  
  
  However $ab=\gamma_m=|\widetilde E(\F_q)|\le 2q$ whereas,  by a previous displayed inequality, $|\Lambda|\ge c_3 q$, a contradiction  for small enough $\delta$, e.g. for $\delta =c_3/8$.\smallskip

  We then reach the crucial conclusion  that  {\it there exist $(t_1,u_1)$ as above with $\gcd(u_1,b_1)=1$}. This corresponds to a point $\tau\in \widetilde E(\F_q)$ such that $\tau=(t_1,u_1)$ in a basis as above, and such that $\tau=[B]x$ for an $x\in\Lambda$. So $\tau$ has exact order $b_1$ and does not lift to $\widetilde W(\F_q)$.\smallskip

  We have choosen $m\in P$ so that $p$ does not divide $\gamma_m=ab$, so the torsion points in $\widetilde E(\F_q)$ may be lifted to $\overline \Q$, and we get a torsion point $\theta\in E[b_1]$ reducing to $\tau$ modulo some place $v'$ of $k_1(E[b_1])$ above $v$, and such that $\theta$  does not lift to $W(k_1(\theta))$:  if this last fact was untrue,  we could reduce modulo $v'$ and obtain a contradiction. (In fact, the restriction $v''$ of $v'$ to $k_1(\theta)$ has residue field $\F_q$ over $v$, because the $m$-th power of the  Frobenius of $v'|v$ fixes   $k_1(\theta)$. Therefore the reduction of $\theta$, i.e. $\tau$, would lift to $\widetilde W(\F_q)$.)

  We can now conjugate over $k_1$ and obtain that all such conjugates $\theta^g$ have the same properties as $\theta$: namely, {\it $\theta^g$ does not lift to $W(k_1(\theta^g))$}.
  
  Now, recall again that we have chosen $m\in P$, so that $\gamma_m$ is not divisible by a fixed large power $l^e$, $e\ge c_2$,  of a prime $l\le c_1$ or $l=l_i$; in particular, $\widetilde E(\F_q)$ does not contain points of  order $l ^e$ for such $l$. As a consequence,  
  we have that, if $m_0$,  and hence $N$,  is large enough, $b_1$ is coprime to such $l$, for otherwise  the $l$-part of $B$ would divide $b$, and $\widetilde E(\F_q)$ would have a  point of order $l^{u+1}$, where $l^u|| B$; but $u\ge N/2l$, so this cannot hold if $N>4c_1c_2\max (l_1,\ldots ,l_{n-1})$. 
  
  But then, by ST, for large enough $c_1$  the set of such conjugates $\theta^g$ over $k_1$ is the whole set of all torsion points of exact order $b_1$. Hence, since $q,B$ were chosen independently of the groups $H$,  we may choose the point $\theta$ independently of $W$ among the $W_H$.\note{Note that this conclusion works for the  points over $\overline\Q$; we couldn't have achieved it over $\F_q$.}

  Taking into account the definition of the $W_H$, and setting $\zeta_n:=\theta$, the proof of   Lemma 1 is thus concluded. \CVD \medskip
  \bigskip

This lemma immediately implies another conclusion:\medskip

\noindent{\bf Lemma 2.}  {\it For $\zeta$ as in Lemma 1, any point $\rho\in Y$  in the fiber $\pi^{-1}(\zeta)$  has degree $\deg\pi>1$ over $k(\zeta)$. Also, there is an automorphism $g \in \Gal(\overline\Q/k(\zeta))$ which does not fix  any   point  in such fiber.}\medskip

\noindent{\it Proof.} For the proof of the first assertion, it suffices to take   into account the  properties of the $Y_H$, as recalled just before the statement of Lemma 1. For the second assertion, just recall Jordan's observation (as in the proof of Theorem 3.1) that any finite  transitive  permutation group not reduced to the identity has an element without fixed points. \CVD
  \medskip 
  
  \noindent{\it End of the proof of Proposition 4.1.}  We let  $m$ be the order of $\zeta$ and $\rho_1,\ldots ,\rho_d$ be the  points in $\pi^{-1}(\zeta)\subset Y$; we put 
  $$
  K:=k(\zeta,\rho_1,\ldots ,\rho_d)
  $$   
  and we let $\tau\in \Gal(K/k(\zeta))$ be the restriction of a $g$  as in the conclusion of Lemma 2, so $\tau$ does not fix any $\rho_i$.  
    
  Now, the strategy will be analogous to that for Proposition  3.1: first, to find primes $l$ such that  $\zeta$ does not lift  to $T$  modulo some place above $l$, and, second,  such that the reduction of $\zeta$ lies in  $\Omega$. This will require some preliminaries. \medskip

  First, letting $\xi=(\xi_1,\ldots ,\xi_n)$ be a generator of $\Omega$, we note that the $\xi_i$ are $\Z$-linearly independent points because $\Omega$ is Zariski dense. Therefore for large $c_1$   the $\xi_i$ are independent also modulo $[m]E(k)$ (actually modulo $[l]E(k)$, for every prime $l>c_1$): a simple proof of this known fact is as in  Proposition 3.1, for the multiplicative case, using N\'eron-Tate heights in place of Weil heights. 
  
  We now pick algebraic points $\eta_1,\ldots ,\eta_n\in  E(\overline k)$ such that
  $[m]\eta_i=\xi_i$ and we put $\eta:=(\eta_1,\ldots ,\eta_n)$.  
   We also put $\zeta=(\zeta_1,\ldots ,\zeta_n)$ and  set $Z:=\Z\zeta_1+\ldots +\Z\zeta_n$, a subgroup of $E[m]$, of exact exponent $m$.  We may find a basis $t,t'$ of $E[m^2]$ such that $Z$ is generated by $mt, mat'$ for a certain divisor $a>0$ of $m$: $m=ab$, say. Note that $Z$ contains $E[b]$ and we have $|Z|=mb$.

   Let us look at the subgroup $H$ of $GL_2(\Z/m^2)$ fixing pointwise $Z$. This subgroup corresponds by ST to the extension $k(E[m^2])/k(\zeta)$. In matrix representation (with respect to the basis $t,t'$), $H$ consists of the $2\times 2$-matrices $I+M$ 
  over $\Z/m^2$, invertible mod $m$ and  such that the first column of $M$ is divisible by $m$ and the second one by $b$. The determinant map taken modulo $m$ gives a homomorphism of $H$ into $(\Z/m)^*$ whose kernel is easily checked to have    cardinality  $m^4a$. Thus any possible intermediate field $F$ between $k(\zeta)$ and $k(E[m^2])$, and whose degree over $k(\zeta)$ is   $<c_1$,  must  be contained in the field corresponding to this kernel. The Galois group of this last field over $k(\zeta)$ is   isomorphic to a subgroup of $(\Z/m)^*$, namely  to the    image $\det H$ (modulo $m$) of $H$ by the determinant homomorphism; in turn, this image is   the group of classes mod $m$ coprime to $m$ and $\equiv 1\pmod b$. If $m'$ is the part of $m$ made up with the primes coprime to $b$ (possibly $m'=1$), then the  subgroup of $\det H$ corresponding to $F$  is  identified with  a subgroup of index $[F:k(\zeta)]$ in $(\Z/m')^*$.   Let now $r$ be a class modulo $ma$, coprime to $m$ and such that $1+br$ is  coprime to $m$ and a primitive root modulo each  prime  dividing $m'$ (recall $\gcd(m',b)=1$). The Chinese Remainder Theorem  delivers such an $r$. Define now the matrix $\mu$  over $\Z/m^2$ by
  $$
  \mu:=I+b\theta,\qquad \theta:=\pmatrix{a\ 0\cr 0\ r}\in {\rm Mat}_2(\Z/ma).\eqno(4.1)
  $$
   As above,  by ST, $\mu$ corresponds to an element of $H\cong \Gal(k(E[m^2])/k(\zeta))$; by our choice of $r$, the group generated by this element has index in $H$ divisible only by primes $>c_1$, and thus cannot fix any field $F$ as above.    \smallskip
   
   Consider now the matrix $\hat\theta:=\pmatrix{r\ 0\cr 0\ a}\in {\rm Mat}_2(\Z/ma)$. It satisfies $\hat\theta+\theta={\rm Tr}(\theta)\cdot I$ and $\theta\hat\theta=\pmatrix{ar\ 0\cr 0\ ar}=\det(\theta)\cdot I$  modulo $ma$. Also, in the said representation of $E[m^2]$ with basis $t,t'$, the image $\hat\theta(E[m])$ is well defined (namely, even if $\hat\theta$ is defined only modulo $ma$ rather than modulo $m^2$) and equals precisely $Z$. Hence there are $z_1,\ldots ,z_n\in E[m]$ such that
   $$
   \hat\theta(z_i)=\zeta_i.\eqno(4.2)
   $$

  Now we recall some facts on elliptic Kummer Theory, relying crucially on ST, for which we refer e.g.  to the paper [Be] and to [L], V (\S 5: `Bashmakov's Theorem'). (In [Be]  the case of prime order is treated; we also mention  papers of Bashmakov, quoted in [Be], [L].)
  
  Since  $\xi_1,\ldots ,\xi_n$ are independent modulo $[m]E(k)$,  we may find, for every choice of points $z_1,\ldots ,z_n$ in $E[m]$, an  automorphism $\sigma\in \Gal(\overline k/k(E[m]))$ such that 
  $$
  \eta_i^\sigma=\eta_i+z_i,\qquad i=1,\ldots ,n.\eqno(4.3)
  $$

  Also, since the Kummer Theory continues to be valid   by replacing $m$ with $m^2$,   the extensions  
  $k(E[m^2])/k(\zeta)$ and $k(\eta)/k(\zeta)$ are linearly disjoint. In fact, the degree of $\eta$ over $k(E[m^2])$ continues to be $m^{2n}$, like the degree  over $k(E[m])$. This degree  is  divisible only by primes $>c_1$, whereas $[K:k(\zeta)]$ divides $\deg\pi !=d!$.  So for large $c_1$ we have  linear disjointness of $k(\eta)/k(\zeta)$ and $K(E[m^2])/k(\zeta)$.    Since, by the above choice and remarks,   the automorphism corresponding to $\mu$    does not fix in particular any nontrivial subextension of $(K\cap k(E[m^2]))/k(\zeta)$,  there exists $g\in\Gal (\overline k/k(\zeta))$ such that
  $$
  g|_{k(E[m^2])}=\mu=1+b\theta\in \Gal(k(E[m^2])/k(\zeta)),\qquad g|_K=\tau|_K,\qquad  \eta^g=\eta^\sigma.\eqno(4.4)
  $$

  Now, we perform a number of choices and deductions:\medskip
    
  1. Let us   choose a large prime $l$ (coprime to $m$ and of good reduction for $E$)   so that its Frobenius with respect to some place $v|l$ of the normal closure $L$ of $K(\eta,E[m^2])/k$ is $g$. This is possible by Chebotarev Theorem. Note in particular that by reduction modulo $v$, $\zeta$ reduces inside $\widetilde E^n(\F_l)$, since $g$ fixes $k(\zeta)$.\smallskip
  
  2. Denote as above by $\widetilde E$ the reduction of $E$ modulo $v$, and let  $\varphi$ be the   Frobenius on $\widetilde E$. Since by n. 1 the point $\zeta$ reduces to a point   defined over $\F_l$, the reduction of the group $Z$ is contained in the kernel of $\varphi -1$, so in particular the kernel of $[b]:\widetilde E\to \widetilde E$ is contained in the kernel of $\varphi-1$. Hence we may write   $\varphi=1+b\psi$, for some endomorphism $\psi$ of $\widetilde E$. (See [Si], Cor. 4.11, p. 77.)
  \smallskip

  3. By n. 1 and (4.4), the Frobenius of $v|l$ acts as $1+b\theta$ on $E[m^2]$, so $b\psi$ and $b\theta$ have the same action on $E[m^2]$, i.e., $b(\psi-\theta)E[m^2]=0$, so in particular $(\psi-\theta)bE[m^2]=0$ and $\psi=\theta$ on $E[m]$.\smallskip

  4. By reduction modulo $v$ of (4.3),  by n. 1 and (4.4), we have on $\widetilde E$, $(\varphi-1)\eta_i=z_i$,  so by  n. 2, $b\psi (\eta_i)=z_i$, (where we have denoted the reduced points by the same letters). \smallskip
  
  5. Let us apply to  the last equations the dual endomorphism $\hat\psi$ of $\psi$, to obtain (recalling $\hat\psi\cdot\psi=\deg(\psi)$)  that $b\deg(\psi)(\eta_i)=\hat\psi(z_i)$, $i=1,\ldots ,n$.  \smallskip
  
  6. Now,  by n. 3, $\psi$ acts as $\theta$  on $E[m]$, and so (by general theory) the same is true of   $\hat\psi$ and $\hat\theta$, so $\hat\psi(z_i)=\zeta_i$ by (4.2), and the equations at n. 5 become
   $b\deg(\psi)\eta\equiv\zeta\pmod v$. Finally, $\deg(\psi)$ is divisible by $a$ (by (4.1) and n. 3), so $b\deg(\psi)=qm$ where $q\in\Z$ and we get 
   $$
   q\xi\equiv \zeta\pmod v.\eqno(4.5)
   $$
  (The fact that $\zeta$ may be represented by some multiple of $\xi$ modulo $v$ is here crucial and seems to be not entirely free of independent interest.) 
  
  Now we can conclude as follows. Let $\gamma=[q]\xi$, $\Omega'=[qm]\Omega=[qm\Z]\xi$. Then, by the congruence (4.5) we have $x\equiv \zeta\pmod v$ for $x\in\gamma+\Omega'$. We suppose that $l$ has been chosen large enough so that $Y$ has good reduction at $v$. Then, if $x$ lifts to a rational point $y\in Y(k)$,  we may reduce modulo $v$ and  obtain that $\zeta$ lifts to $\widetilde Y(\F_l)$, so the reduction of some $\rho_i$ lies in $\F_l$. This is however impossible because by (4.4) and n. 1 above the Frobenius of $v|l$ acts as $g$ on $K$ and hence it moves each $\rho_i$.
  
  This proves the case $h=1$ of Proposition 4.1 (with $C=\gamma+\Omega'$) and, as remarked in the opening argument, this suffices for the general case. \CVD
   \bigskip

   \noindent{\bf Remarks.}  
  (i)  
  One may argue similarly with  CM curves (see [L]), which we leave  to the interested reader, and with non cyclic finitely generated $\Omega$, with more complicated, but conceptually similar, arguments.  (The cyclic case is the most basic here, also because the main purpose in Hilbert Irreducibility is to find   `good' elements in `small' sets.) 
  
    (ii) As already remarked,  the method of proof in principle applies to more general abelian varieties, provided one has the suitable  Kummer Theory and torsion-Galois action at disposal. This seems not to be yet available in the most general case, but probably one can deal with other special cases.  All of this goes in the direction of  the problem stated at p. 53, \S 5.4 of [Se2]. 
    
    (iii) On weakening assumption (PB) one may obtain corresponding versions of Theorems 1,2 in which `irreducible' is replaced by `irrational'. More precisely, denote  by $X$ either $\G_{\rm m}^r\times \G_{\rm a}$ or $E^n$. We have:\medskip
    
    \noindent{\bf Theorem 4.} {\it For $i=1,\ldots ,h$, let $\pi_i:Y_i\to X$ be a  cover,   not birationally equivalent to an isogeny. Let $\Omega\subset X(k)$ be a Zariski dense cyclic subgroup. Then there exists a coset $C$ of finite index in $\Omega$  and disjoint from $\cup_{i=1}^h\pi_i(Y_i(k))$.}\medskip

   This is easily deduced from Theorem 1 or 2 (depending on $X$). In fact, by Proposition 2.1, each map $\pi_i$ factors as $\lambda_i\circ\rho_i$ where $\lambda_i$ is an isogeny and $\rho_i$ satisfies (PB).  Since none of the covers is birationally equivalent to an isogeny, each $\lambda_i$ has degree $>1$.  The isogeny $\lambda_i$ is a factor of a multiplication map $[m_i]=\lambda_i\circ \hat\lambda_i:X\to X$, so, replacing $Y_i$ by its pullback $\hat\lambda_i^*Y_i$, with $\pi_i$ replaced by  the natural map to $X$ induced by $\rho_i$,  we may suppose that $Y_i$ satisfies (PB). Now, by Theorems 1,2   we deduce the existence of the coset $C$ such that each lifting of a point of $C$ by $\rho_i$ has degree $>1$, proving the claim. (The liftings by the $[m]_i$ do not matter, since they occur over a fixed number field, which can be supposed to be $k$.) \bigskip

\centerline{\sc \S 5. Theorem 3 and an elliptic analogue.}\medskip

We start by proving Theorem 3.
 We recall that  $k^c$ denotes the   extension of $k$ obtained by adding to it all roots of unity. \medskip

\noindent{\it Proof of Theorem 3.}   We first easily reduce to the case $\kappa=\overline\Q$ by means of a specialization argument. 
First, up to birationality we may assume that $\pi$ is finite. Let $\overline\Q\subset\kappa_0\subset\kappa$ be a field of definition for $Y$, $ \kappa_0$ finitely generated over $\overline\Q$. We may view $Y$ as a finite cover, defined over $\overline \Q$,  of  $V\times \G_{\rm m}^n$, where $V$ is an affine variety with function field  $\kappa_0$.  For  $\xi\in V(\overline\Q)$ we have a specialized cover $Y_\xi\to\G_{\rm m}^n$.  Since $Y$ is irreducible over $\kappa$, which is algebraically closed, it is well known that $Y_\xi$ is irreducible over $\overline \Q$ for all $\xi$ in a Zariski-open $U\subset V$ (thinking of $Y$ as defined by a polynomial, this follows e.g. from [Sch], Thm. 32 and Cor. 2, pp. 201-202).  By the same argument applied to the pullback cover $[d]Y$, we may shrink $U$ and suppose that $Y_\xi$ satisfies (PB) and that the degree $[Y_\xi :\G_{\rm m}^n]$ equals $d=\deg\pi$. Pick now $\xi\in U$, let $k$ be a number field of definition for $Y_\xi$ and apply the conclusion of Theorem 3 to $Y_\xi$, obtaining a set $\E$ as therein. We contend that this set works for $Y$ as well. In fact, let $G$ be a connected algebraic subgroup of $\G_{\rm m}^n$, $G\not\subset\E$, and $\theta$ a torsion point.  Suppose that $\pi^{-1}(\theta G)$ is reducible over $\kappa$.  Then, since $\pi$ is finite,  there is an affine variety $V_1$ with a finite map $\rho:V_1\to V$ such that the pullback $\rho^* Y$, as  a cover of $V_1\times G$, is reducible, say as a union $Z\cup W$ of two covers $Z,W$ of degree $<\deg\pi$. Now, since $\rho$ is finite, there is a $\xi_1\in V_1(\overline\Q)$ with $\rho(\xi_1)=\xi$. It follows that $Y_\xi\cap\pi^{-1}(\theta G)$ is  the  union $Z_{\xi_1}\cup W_{\xi_1}$ and thus reducible  (the fact that  $Y_\xi\to\G_{\rm m}^n$
 has degree $\deg\pi$ ensures that $Z_{\xi_1}\cup W_{\xi_1}$ is a nontrivial decomposition), a  contradiction.
\medskip
Let us now prove the theorem in the crucial case $\kappa=\overline\Q$. Let $k$ be a number field of definition for $Y$ and $\pi$, and let us apply, as we may,  Theorem 2.1 to our cover $Y\to \G_{\rm m}^n$, obtaining a finite union $\E_1$ of torsion cosets as therein. By applying that conclusion to torsion points $\zeta\in \theta G\setminus \E_1$ and   recalling that torsion points are Zariski-dense in $G$, we obtain that if $\theta G\not\subset\E_1$,  then $\pi^{-1}(\theta G)$ is irreducible over $k^c$ (for otherwise $\pi^{-1}(\zeta)$ would be {\it a fortiori} reducible over $k^c$ for the Zariski-dense set of torsion points $\zeta\in \theta G\setminus \E_1$).

The point is now to go from $k^c$ to $\bar k$, and for this  we consider the cover $W:=Y\times Y\to \G_{\rm m}^{2n}\cong \G_{\rm m}^n\times \G_{\rm m}^n$, by the map $\pi_2:=\pi\times \pi$ of degree $d^2$ where $d:=\deg\pi$. Since $Y$ satisfies (PB), the same is true of $W$, as a cover of $\G_{\rm m}^{2n}$.  Hence by Theorem 2.1 applied this time to $W,\pi_2$ we deduce that there is a finite union $\E_2$ of proper torsion cosets of $\G_{\rm m}^{2n}$ such that for $\zeta_1\times\zeta_2$ a torsion point in $\G_{\rm m}^{2n}\setminus\E_2$ the fiber $\pi_2^{-1}(\zeta_1\times\zeta_2)$ is $k^c$-irreducible. 

Denote $Z:=\pi^{-1}(\theta G)$  and suppose that $Z$ is reducible over $\bar k$. If $\theta G\not\subset \E_1$, we have observed that $Z$ is irreducible over $k^c$ and then the function field extension $k^c(Z)/k^c(G)$ contains a nontrivial finite `constant' extension $L/k^c$. But then  $Z\times Z$ is reducible over $k^c$, and in fact each $k^c$-component $Z_2$ satisfies $[k^c(Z_2):k^c(G\times G)]\le [k^c(Z):k^c(G)]^2/[L:k^c]=d^2/[L:k^c]$.\note{Note that we cannot directly work over $L$, which {\it a priori}  might depend on $G$. See next theorem for an alternative argument.}   Hence the fiber in  $Z\times Z$ above  a torsion point $\zeta_1\times\zeta_2\in \theta G\times \theta G$  has at least $[L:k^c]$ components  irreducible over $k^c$. We conclude that $\theta G\times \theta G$ is contained in $\E_2$. 

Thus if $\pi^{-1}(\theta G)$ is reducible, we obtain  that either $\theta G\subset \E_1$ or $\theta G\times \theta G\subset \E_2$. 
From this we easily deduce that    $G$ is anyway contained in a certain finite union $\E$ of proper connected algebraic subgroups of $\G_{\rm m}^n$, concluding the argument. \CVD \medskip

We end this section with    a result  in the elliptic context, similar but weaker (in that the subgroup is  restricted to a special shape)  than  Theorem  3;  in this situation   we lack Theorem 2.1, so we cannot argue as above.  We let $E/k$ be an elliptic curve without CM, $r\in\N$   and put $A=E^r$. \medskip

\noindent{\bf Theorem 5.1.} {\it Let $\pi:Y\to A^n$ be a cover satisfying (PB). Then there are integers $a_1,\ldots ,a_n\neq 0$ such that the restriction of $Y$ above the subgroup  $B=\{([a_1]x,\ldots ,[a_n]x):x\in  A\}$  is irreducible. Moreover, we can choose the $a_i$ so that $B$ is not contained in any prescribed finite union $\E$ of proper torsion cosets in $A^n$.}\medskip

\noindent{\it Proof.} 
For $s=1,\ldots ,n$, we prove inductively on $n,s$ that there are integers $a_1,\ldots ,a_s\neq 0$ such that $Y$ is irreducible above $\{([a_1]x,\ldots ,[a_s]x):x\in  A\}\times A^{n-s}$, so that moreover this is not a subset of  $\E$. On taking a pullback by $[\gcd(a_i)]|_{A^s}\times{\rm Id}|_{A^{n-s}}$ we may assume that $\gcd(a_i)=1$, in which case this variety is isomorphic to $A^{n-s+1}$. The  assertion is trivial for $s=1$, any $n$,  and  by induction on $n$  we reduce to  $s=2$. Suppose  that for generic $x\in A^{n-1}$, $Y$ becomes reducible above $\{x\}\times A$. Then the product $Y\times_{A^n} V$ is reducible   in $t>1$ components, where $V=C\times A$ and $C$ is a suitable finite normal  cover of $A^{n-1}$ of degree $t$. If $C\to A^{n-1}$  is unramified, then the cover  is birationally equivalent to an isogeny of   abelian varieties; but the reducibility then    violates  our assumptions on $Y$. Therefore it  is branched, say at  $x_0\in A^{n-1}$. But   two distinct components of $Y\times_{A^n}V$ merge above any point in  the branch locus of $V\to A^n$, so the branch locus of $\pi$ contains $\{x_0\}\times A$. 
   By an automorphism of $A^n$ induced by a unimodular $n\times n$ matrix of integers (which does not affect the result) we may assume this is not the case, so $Y_x:=\pi^{-1}(\{x\}\times A)$ is  generically irreducible,  and so it is  irreducible except for a proper closed subset of $x\in A^{n-1}$.\note{This argument also appears in the proof of Proposition 4.1 above.}   
Arguing similarly  with a pullback $[d]^*Y$  in place of $Y$, we conclude that $Y_x\to A$ remains generically irreducible under pullback by $[d]$ and so (by Prop. 2.1) satisfies (PB).  

Put now $B_m:=\{(x_1,\ldots ,x_n)\in A^n: x_n=[m]x_1\}$ and $Y(m)=\pi^{-1}(B_m)$. Suppose that $Y(m)$ is reducible over $\bar k$, but irreducible over $k$. Then its function field contains a nontrivial finite extension $L/k$, necessarily of degree $\le d$. 

 On the other hand, if $B_m$ is not contained in the branch locus of $\pi$, we contend that the discriminant of $L/\Q$ can be divisible only by primes in a finite set independent of $m$ (but dependent only on $Y,\pi ,k$). In fact, note first that  $Y(m)\to B_m$ is unramified above a generic point and thus its reduction at  a  place $v$ of $k$ may be generically ramified only if  
 the reduction of $B_m$ is contained in the branch locus of  the reduction of $\pi  $. Now,  by reducing modulo $v$ a nontrivial algebraic equation valid on the branch locus of $\pi$, this  yields a fixed  algebraic relation modulo $v$  between $x_1,[m]x_1$. So, since $\deg [m]\to\infty$, for large $m$ this relation must be trivial modulo $v$, and hence $v$ must lie  in a finite set independent of the integers $m$ in question. 
 
 Since   $[L:\Q]\le d[k:\Q]$, we conclude that the discriminant of $L$ is also bounded (see  [Se5], p. 67, Remarque), and thus by a well-known result of Hermite   $L$ has only finitely many possibilities independent of $m$. Let then $k_1$ be the number field generated by all such possible fields $L$. We have proved that in any case either (i) $Y(m)$ is irreducible, or (ii) $B_m$ is contained in the branch locus of $\pi$, or (iii) $Y(m)$ is reducible over a number field $k_1$ independent of $m$. 
 
 Finally, choose $x_1\in A(k_1)$ so that $\Z x_1$ is  Zariski-dense in $A$  (it exists if we enlarge $k_1$). Replacing $x_1$ with a multiple and choosing other points $x_2,\ldots ,x_{n-1}\in A(k_1)$ and setting $z_0:=(x_1,\ldots ,x_{n-1})\in A^{n-1}$, we have seen that we may assume that the variety $Y_{z_0}$ and its pull-back by $[d]$ are  irreducible. Then,  by Theorem 2  applied to $Y_{z_0}\to A$  there exist infinitely many $m\in\N$ such that the fiber above $[m]x_1$ in $Y_{z_0}$, namely, $\pi^{-1}((z_0,[m]x_1))$ is $k_1$-irreducible. But $(z_0,[m]x_1)\in B_m$, so $Y(m)$ must be itself $k_1$-irreducible. Hence alternative (iii) is not verified, and the same holds for (ii) if $m$ is large. Therefore $Y(m)$ is irreducible for infinitely many $m$. On the other hand, it is easily checked that $B_m$  can be contained in $\E$ only for finitely many $m$,  which proves the result.\CVD
 \medskip

\noindent{\bf Remark.}  We could also have applied Thm. 2  to a product $Y_{z_0}\times Y_{z_1}$, with the argument for Thm. 3 in place of exploiting ramification. Also, for $r=1$ one can use   the function field version of the Mordell Conjecture in place of Thm. 2.
\bigskip

 \bigskip
\centerline{\title References}\medskip

\item{[Be]} - D. Bertrand, {\it Kummer Theory on the product of an elliptic curve by the multiplicative group}, Glasgow Math. J., {\bf 22} (1981), 83-88.\smallskip

\item{[BG]} - E. Bombieri, W. Gubler, Heights in Diophantine Geometry, {\it Cambridge Univ. Press}, 2006\smallskip

\item{[BM]} - W.D. Brownawell,  D. Masser, {\it Vanishing sums in function fields}, Math. Proc. Camb. Phil. Soc., {\bf  100}  (1986), 427-434.\smallskip

\item{[CS]} - J.-L. Colliot-Th\'el\`ene, J.-J. Sansuc, {\it Principal Homogeneous Spaces under Flasque Tori: Applications}, J. Algebra, {\bf 106} (1987), 148-205.\smallskip

\item{[C]} - P. Corvaja, {\it Rational fixed points for linear group actions}, Ann. Scuola Norm. Sup. Pisa Cl. Sci., Ser. V, {\bf 6} (2007),  561-597.\smallskip

\item{[CZ]} . P. Corvaja, U. Zannier, {\it Some new applications of the subspace theorem}, Compositio Math. {\bf 131} (2002), no. 3, 319-340.\smallskip

\item{[CZ2]} . P. Corvaja, U. Zannier, {\it On the Integral Points on Certain Surfaces}, Intern. Math. Res. Notices, 2006, 1-20.\smallskip

\item{[CZ2]} . P. Corvaja, U. Zannier, {\it  Some cases of Vojta's conjecture on integral points over function fields}, J. Alg. Geometry, {\bf 17} (2008), 295-333.\smallskip

\item{[D]} - P. D\`ebes, {\it On the irreducibility of the polynomials $P(t\sp m,Y)$},  J. Number Theory {\bf 42} (1992), no. 2, 141-157.\smallskip

\item{[DZ]} - R. Dvornicich, U. Zannier, {\it Cyclotomic diophantine problems (Hilbert Irreducibility and invariant sets  for polynomial maps)}, Duke Math. J., {\bf 139} (2007).\smallskip

\item{[FZ]} - A. Ferretti, U. Zannier, {\it Equations in the Hadamard ring of rational functions}, Ann. Scuola Norm.  Sup. Pisa, Cl. Sci., Ser. V, {\bf 6} (2007), 457-475.\smallskip

\item{[Fr]} - M. Fried, {\it On Hilbert's Irreducibility Theorem}, J. Number Theory, {\bf 6} (1974), 211-231.\smallskip

\item{[FJ]} - M. Fried, M. Jarden, Field arithmetic, {\it Springer-Verlag}.\smallskip

\item{[K]} - S. Kleiman, {\it The transversality of a general translate}, Compositio Math., {\bf 28} (1974), 287-297.\smallskip

  \item{[I]} - A.E. Ingham, The Distribution of Prime Numbers, {\it Cambridge Univ. Press}, 1932, Reprinted 1992.\smallskip

\item{[L]} - S. Lang, Elliptic curves, Diophantine Analysis, {\it Springer Verlag}, 1978.\smallskip

\item{[L2]} - S. Lang, Number Theory III, Diophantine Geometry, Enc. of Mathematical Sciences, Vol. 60, {\it Springer Verlag}, 1991.\smallskip

\item{[Sch]} - A. Schinzel, Polynomials with special regard to reducibility, {\it Cambridge Univ. Press}, 2000.\smallskip

\item{[Sch2]} - A. Schinzel, {\it On Hilbert's Irreducibility Theorem}, Ann. Polonici Math.,  {\bf XVI} (1965), 333-340.\smallskip

\item{[Sch3]} - A. Schinzel, {\it An analogue of  Hilbert's irreducibility theorem}, Number Theory, W. de Gruyter, Berlin (1990), 509-514.\smallskip

\item{[Se1]} - J-P. Serre, Lectures on the Mordell-Weil Theorem,  $2$-nd ed., {\it Vieweg}, 1990.\smallskip

\item{[Se2]} - J-P. Serre, Topics in Galois Theory, {\it Jones and Bartlett}, Boston, 1992.\smallskip

\item{[Se3]} - J-P. Serre, {\it On a theorem of Jordan}, Bull. A.M.S., {\bf 40} (2003), 429-440.\smallskip

\item{[Se4]} - J-P. Serre,  {\it Propri\'et\'es galoisiennes des points d'ordre fini des courbes elliptiques},  Invent. Math. {\bf 15} (1972), no. 4, 259-331.\smallskip

\item{[Se5]} - J-P. Serre, Corps locaux, Tr. \'ed., {\it Hermann}, Paris, 1968.\smallskip

\item{[Si]} - J. Silverman, The Arithmetic of Elliptic Curves, {\it Springer Verlag}, GTM {\bf 106}, 1986.\smallskip

\item{[V]} - H. Volklein, Groups as Galois Groups, {\it Cambridge Univ. Press}, 1996.\smallskip

\item{[Z]} - U. Zannier, {\it A proof of Pisot $d$-th root conjecture},  Annals Math. {\bf 2} (2000), 375-383.\smallskip

\bigskip

\vfill 

Umberto Zannier

Scuola Normale Superiore

Piazza dei Cavalieri, 7

56126 Pisa - ITALY \hfill e-mail:  {\it u.zannier@sns.it}

\end